\documentclass[12pt]{amsart}
\font\emailfont=cmtt10

\usepackage{amsmath,amsthm,amsfonts,amscd,flafter,epsf,latexsym}

\newtheorem{theorem}{Theorem}[section]
\newtheorem{prop}[theorem]{Proposition}
\newtheorem{cor}[theorem]{Corollary}
\newtheorem{lemma}[theorem]{Lemma}

\newtheorem{defn}[theorem]{Definition}

\newtheorem{remark}[theorem]{Remark}

\def\endproof{\relax\ifmmode\expandafter\endproofmath\else
  \unskip\nobreak\hfil\penalty50\hskip.75em\hbox{}\nobreak\hfil\bull
  {\parfillskip=0pt \finalhyphendemerits=0 \bigbreak}\fi}
\def\endproofmath$${\eqno\bull$$\bigbreak}
\def\bull{\vbox{\hrule\hbox{\vrule\kern3pt\vbox{\kern6pt}\kern3pt\vrule}\hrule}}

\newcommand{\R}{\mathbb{R}}

\newcommand{\C}{\mathbb{C}}

\newcommand{\Z}{\mathbb{Z}}

\newcommand{\Projectivize}{\mathbb{P}}
\newcommand{\OneHalf}{\frac{1}{2}}

\newcommand{\CP}[1]{{\mathbb{CP}}^{#1}}
\newcommand{\CPbar}{{\overline{\mathbb{CP}}}^2}
\newcommand{\Zmod}[1]{\Z/{#1}\Z}

\newcommand{\Tr}{\mathrm{Tr}}
\newcommand{\Ker}{\mathrm{Ker}}
\newcommand{\CoKer}{\mathrm{Coker}}

\newcommand{\cm}{\cdot}
\newcommand\Sections{\mbox{$\Gamma$}}

\newcommand{\nbd}[1]{\mathrm{nd}({#1})}

\newcommand{\ModSW}{\ModSWfour}
\newcommand{\ModSWirr}{\ModSW^{\irr}}
\newcommand{\ModSWIrr}{\ModSWirr}

\newcommand{\ModSWthree}{\mathcal{N}}

\newcommand{\ModSWfour}{\mathcal{M}}
\newcommand{\ModFlow}{\ModSWfour}
\newcommand{\UnparamModFlow}{\widehat{\ModFlow}}
\newcommand{\UnparModFlow}{\UnparamModFlow}

\newcommand{\SW}{SW}

\newcommand{\Sobol}[2]{L^{#1}_{#2}}

\newcommand{\Canon}[1]{{{K}_{#1}}}

\newcommand{\Dirac}{\mbox{$\not\!\!D$}}

\newcommand{\Conns}{\mathcal A}

\newcommand{\SWConfig}{\mbox{${\mathcal B}$}}
\newcommand{\SWPreConfig}{\mbox{${\mathcal C}$}}
\newcommand{\SWConfigIrr}{\SWConfig^{\irr}}
\newcommand{\SWPreConfigIrr}{\SWPreConfig^{\irr}}
\newcommand{\irr}{\mathrm{*}}
\newcommand{\red}{\mathrm{red}}

\newcommand{\Map}{\mathrm{Map}}

\newcommand{\SpinC}{{\mathrm{Spin}}_{\C}}

\newcommand{\Proj}{\Pi}

\newcommand{\dd}[2]{\frac{d{#1}}{d{#2}}}
\newcommand{\DD}[1]{\frac{\partial}{\partial{#1}}}
\newcommand{\DDs}{\DD{s}}
\newcommand{\DDt}{\DD{t}}

\newcommand{\Deriv}{D}
\newcommand{\Vol}{\mathrm{Vol}}

\newcommand{\goesto}{\mapsto}

\newcommand{\DBar}{\overline{\partial}}
\newcommand{\Del}{\partial}

\newcommand{\CSD}{\mbox{${\rm CSD}$}}

\newcommand\edim{\text{\rm{e-dim}}}

\newcommand\Wedge{\Lambda}

\newcommand\loc{\mathrm{loc}}

\newcommand\Hom{\mathrm{Hom}}

\newcommand\abuts\Rightarrow
\newcommand\Sym{\mathrm{Sym}}

\newcommand\Jac{\mathcal{J}}
\newcommand\Alg{\mathbb{A}}

\newcommand\Restrict{\rho}

\newcommand{\UnparModFlowBased}{\UnparModFlow^0}
\newcommand\ModSWBased{\ModSW^0}

\newcommand\Hol{\mathrm{Hol}}
\newcounter{bean}

\newcommand\NormBun{L}
\newcommand\NormBundle{\NormBun}

\newcommand\ESigma{E_0}

\newcommand\RSurf{R}

\newcommand{\PD}{\mathrm{PD}}

\newcommand{\ModSp}{\ModSWfour}
\newcommand{\ModSWThree}{\ModSWthree}

\newcommand{\ModRed}{\ModSp^{\red}}
\newcommand{\ModIrr}{\ModSp^{\irr}}

\newcommand\spinc{\mathfrak s}
\newcommand\spinct{\mathfrak t}

\newcommand\Xtrunc{X_0}
\newcommand\TConn{\nabla}

\newcommand{\SWirr}{\SW^{irr}}
\newcommand{\SWred}{\SW^{red}}
\newcommand{\CritMan}{C}

\newcommand{\Cinfty}{C^{\infty}}
\newcommand{\CinftyLoc}{\Cinfty_{\loc}}
\newcommand{\mult}{e}
\newcommand{\Euler}{\mathbf{e}}

\newcommand{\CompactModSW}{\overline\ModSW}
\newcommand{\CompactRestrict}{\overline\rho}
\newcommand{\Include}{i}
\newcommand{\CompactInclude}{\overline{i}}

\newcommand{\BaseFib}{\mathcal L}
\newcommand{\SWTube}{{\widehat\SW}_{(\CritMan,\Jac)}}
\newcommand{\Id}{\mathrm{Id}}
\newcommand{\NSig}{N}

\newcommand{\Divisors}{V}
\newcommand{\QuotMap}{q}

\newcommand{\BigSpinc}{-n-2g+2 \leq \langle c_1(\spinc),[\Sigma]\rangle\leq -n}
\newcommand{\SmallSpinc}{-n<\langle c_1(\spinc),[\Sigma]\rangle\leq  n}
\newcommand{\OBundle}{\Xi}

\newcommand{\Chamber}{\mathcal K_0}
\newcommand{\Chambers}{\mathcal K}

\title{Higher Type Adjunction Inequalities in Seiberg-Witten Theory}
\author[Peter Ozsv{\'a}th]{Peter Ozsv\'ath} \thanks{The first author was partially supported by NSF grant number DMS 9304580}
\address{School of
Mathematics, Institute for Advanced Study, New Jersey 08540 \newline
\indent{\emailfont{petero@math.ias.edu}}}

\author[Zolt{\'a}n Szab{\'o}]{Zolt{\'a}n Szab{\'o}} 
\thanks{The second author was partially 
supported by NSF grant number DMS 970435 and a
Sloan Fellowship} 
\address{Department of
Mathematics,  Princeton University, New Jersey 08540 \newline
\indent{\emailfont{szabo@math.princeton.edu}}}

\begin{document}

\begin{abstract} 
In this paper, we derive new adjunction inequalities for embedded
surfaces with non-negative self-intersection number in
four-manifolds. These formulas are proved by using relations between
Seiberg-Witten invariants which are induced from embedded surfaces.
To prove these relations, we develop the relevant parts of a
Floer theory for four-manifolds which bound circle-bundles over
Riemann surfaces.
\end{abstract}

\maketitle
\section{Introduction}
\label{sec:Introduction}

In this paper, we prove certain adjunction inequalities, which give
relations between the Seiberg-Witten invariants of a four-manifold $X$
and the genus of embedded surfaces in $X$. These results are
generalizations of results from ~\cite{KMthom},
\cite{MSzT}, \cite{SympThom}, see also~\cite{KMStructureTheorem}.

The investigations center on a construction of an appropriate
Seiberg-Witten-Floer functor for manifolds which bound circle bundles
$Y$ over Riemann surfaces (with sufficiently large Euler number),
which relies on the calculations of~\cite{MOY}. Special cases of this
theory were studied in~\cite{SympThom}, where the authors used similar
techniques to prove the symplectic Thom conjecture. That problem
requires an analysis of those $\SpinC$ structures over $Y$ for which
the Seiberg-Witten moduli space contains only reducible solutions,
which simplifies the corresponding Floer homology. In this paper, we
work out the theory in the other, more complicated cases. We will give
more applications of these techniques in~\cite{OSzThree}.

Before stating the results, we set up some notation.  Let $X$ be a
closed, connected, smooth four-manifold equipped with an
orientation for which $b_2^+(X)>0$ (where $b_2^+(X)$ is the dimension
of a maximal positive-definite linear subspace $H^+(X;\R)$ of the intersection
pairing on $H^2(X;\R)$) and an orientation for $H^1(X;\R)\oplus
H^+(X;\R)$. Given such a four-manifold, together with a $\SpinC$
structure $\spinc$, the Seiberg-Witten invariants (see~\cite{Witten},
\cite{Morgan},
\cite{Salamon}) form an integer-valued function $$\SW_{X,\spinc}\colon
\Alg(X)\longrightarrow
\Z,$$
where $\Alg(X)$ denotes the graded algebra obtained by tensoring
the exterior algebra on $H_1(X)$ (graded so that $H_1(X)$ has grading
one) with the polynomial algebra $\Z[U]$ on a single two-dimensional
generator. The invariants are constructed via intersection theory on
the moduli space $\ModSW_X(\spinc)$ 
of solutions $(A,\Phi)$ modulo gauge to the Seiberg-Witten equations in
$\spinc$:
\begin{eqnarray}
\rho(F_A^+) &=& i\{\Phi,\Phi\}_0 - \rho(i\eta) \label{eq:CurvatureEquation}
\\
\Dirac_A \Phi &=& 0,
\label{eq:HarmonicSpinor}
\end{eqnarray}
where $\Phi$ is a section of $W^+$, $A$ is a spin-connection in the
spinor bundle $W^+$ of $\spinc$, $\Dirac_A$ denotes the associated
Dirac operator, $\rho$ denotes Clifford multiplication, $\eta$ is some
fixed self-dual two-form, and $\{\Phi,\Phi\}_0$ is the usual quadratic
map (see~\cite{Witten}).  Note that the invariants are zero on
homogeneous elements whose degree is not $d(\spinc)$, where
$$d(\spinc)=\frac{c_1(\spinc)^2 - (2\chi(X) +3\sigma(X))}{4}$$ denotes
the formal dimension of the moduli space $\ModSW_X(\spinc)$.  When
$b_2^+(X)>1$, $\SW_{X,\spinc}$ is a diffeomorphism invariant of the
four-manifold; when $b_2^+(X)=1$, the invariants depend on a chamber
structure (see~\cite{Morgan},
\cite{SympThom}).  There are two distinguished chambers
corresponding to the two components of $\Chambers(X)=\{\omega\in
H^2(X;\R)-0~|~\omega^2 \geq 0 \}$. Given a component $\Chamber$ of
$\Chambers(X)$, the corresponding invariant (still denoted $\SW_{X,\spinc}$)
is calculated using the
moduli space of solutions to the Seiberg-Witten equations perturbed by
any generic self-dual two-form $\eta$, provided that the sign of
$-2\pi c_1(\spinc)\cm\omega_g + \int_X \eta\wedge\omega_g$ agrees with
the sign of $\gamma\cm\omega_g$, where $\gamma$ is any class in
$\Chamber$, and $\omega_g\neq 0$ is a harmonic (with respect to the
metric $g$), self-dual two-form over $X$.  Note that $\SW_{X,\spinc}$
is a diffeomorphism invariant of $X$ (and the component $\Chamber$).

Those $\SpinC$ structures $\spinc$ for which the invariant
$\SW_{X,\spinc}$ is non-trivial are called {\em basic classes}.

Our results are easiest to state when $b_1(X)=0$, where we have the
following.

\begin{theorem}
\label{thm:AdjunctionInequalityBOneZero}
Let $X$ be a smooth, closed, connected, oriented four-manifold with $b_2^+(X)>0$
and $b_1(X)=0$, 
and let
$\Sigma\subset X$ be a smoothly-embedded surface 
with genus $g(\Sigma)>0$ representing a
non-torsion homology class with  self-intersection
number $[\Sigma]\cm[\Sigma]\geq 0$. 
If $b_2^+(X)>1$, then we have the following adjunction inequality
$$
|\langle c_1(\spinc),[\Sigma]\rangle|  + [\Sigma]\cm [\Sigma] +
2d(\spinc) \leq 2g(\Sigma)-2,
$$
for each basic class $\spinc\in \SpinC(X)$. Furthermore, when
$b_2^+(X)=1$, for each basic class $\spinc$ of $X$ for the component
of $\Chambers(X)$ which contains $\PD[\Sigma]$
with
$$-\langle c_1(\spinc),[\Sigma]\rangle+[\Sigma]\cm[\Sigma]\geq 0,$$
we have an inequality
$$
-\langle c_1(\spinc),[\Sigma]\rangle  + [\Sigma]\cm [\Sigma] +
2d(\spinc) \leq 2g(\Sigma)-2.
$$
\end{theorem}

\begin{remark}
The above theorem should be seen as a refinement of the adjunction
inequality proved by Kronheimer-Mrowka and Morgan-Szab{\'o}-Taubes
(see~\cite{KMthom}, \cite{MSzT}, \cite{DonaldsonBulletin}). 
Analogous results for immersed spheres were obtained by Fintushel and
Stern, see~\cite{FS}. 
\end{remark}

In fact, Theorem~\ref{thm:AdjunctionInequalityBOneZero} follows from a more
general version.  To state this, note first that an inclusion $i\colon
\Sigma\longrightarrow X$ induces a map
$$i_*\colon \Alg(\Sigma)\longrightarrow \Alg(X).$$

\begin{theorem}
\label{thm:StrongAdjunctionInequality}
Let $X$ be a smooth, closed, connected, oriented four-manifold with
$b_2^+(X)>0$.  Let $\Sigma\subset X$ be a surface with genus
$g(\Sigma)>0$ representing a non-torsion homology class with
self-intersection number $[\Sigma]\cm[\Sigma]\geq 0$. Let $\ell$ be an
integer so that there is a symplectic basis $\{A_j,B_j\}_{j=1}^g$ for
$H_1(\Sigma)$ so that $i_*(A_j)=0$ in $H_1(X;\R)$ for $i=1,...,\ell$.
Let $a\in\Alg(X)$
and $b\in\Alg(\Sigma)$ be an element of degree $d(b)\leq \ell$.  
If $b_2^+(X)>1$ then for each $\SpinC$ structure
$\spinc$ so that $\SW_{X,\spinc}(a\cm i_*(b))$ is non-zero, we have $$
|\langle c_1(\spinc),[\Sigma]\rangle| + [\Sigma]\cm [\Sigma] + 2d(b)
\leq 2g(\Sigma)-2.  $$ Furthermore, when $b_2^+(X)=1$ then for each
$\SpinC$ structure $\spinc$ of $X$ with $$-\langle
c_1(\spinc),[\Sigma]\rangle+[\Sigma]\cm[\Sigma]\geq 0,$$ for which
$\SW_{X,\spinc}(a\cm i_*(b))$ is non-zero, when calculated in the
component of $\Chambers(X)$ containing $\PD[\Sigma]$, we have an
inequality
\begin{equation}
\label{eq:StrongAdjunctionInequality}
-\langle c_1(\spinc),[\Sigma]\rangle  + [\Sigma]\cm [\Sigma] +
2d(b) \leq 2g(\Sigma)-2.
\end{equation}
\end{theorem}

The Adjunction Inequality~\eqref{eq:StrongAdjunctionInequality} does not
hold without homological restrictions on $X$, as we can see by looking
at the ruled surface $X=S^2\times\Sigma$. 
In general, one can obtain only a weaker inequality (losing the
factor of $2$ on the dimension $d(b)$), as follows.

\begin{theorem}
\label{thm:GeneralAdjunctionInequality}
Let $X$ be a smooth, closed, connected, oriented four-manifold with
$b_2^+(X)>1$.  Let $\Sigma\subset X$ be a surface with genus
$g(\Sigma)>0$ representing a non-torsion homology class with
self-intersection number $[\Sigma]\cm[\Sigma]\geq 0$. Let
$a\in\Alg(X)$ and $b\in\Alg(\Sigma)$.  If $b_2^+(X)>1$ and if
$\SW_{X,\spinc}(a\cm i_*(b))$ is non-zero for some $b\in\Alg(\Sigma)$
of degree $d(b)$, then we have $$ |\langle
c_1(\spinc),[\Sigma]\rangle| + [\Sigma]\cm [\Sigma] + d(b) \leq
2g(\Sigma)-2.  $$ If $b_2^+(X)=1$ and $\spinc$ is a $\SpinC$ structure
with $$-\langle c_1(\spinc),[\Sigma]\rangle+[\Sigma]\cm[\Sigma]\geq
0,$$ for which $\SW_{X,\spinc}(a\cm i_*(b))$ is non-zero, when
calculated in the component of $\Chambers(X)$ containing
$\PD[\Sigma]$, then we have
\begin{equation}
\label{eq:GeneralAdjunctionInequality}
-\langle c_1(\spinc),[\Sigma]\rangle + [\Sigma]\cm [\Sigma] +
d(b) \leq 2g(\Sigma)-2.
\end{equation}
\end{theorem}

\begin{remark}
Adjunction inequalities for surfaces of positive square in Donaldson's
theory were first obtained in the influential paper of Kronheimer and
Mrowka (see~\cite{KMStructureTheorem}). These inequalities were strengthened
under similar, but more restrictive, hypotheses in their
preprint~\cite{KMunpublished}; see also~\cite{FSdonaldson}. The
conjectured relationship between the Donaldson and Seiberg-Witten
invariants gives a correspondence between the adjunction inequalities
arising in these two theories. For more on this correspondence,
see~\cite{SeibergWitten}, \cite{Witten}, \cite{WittenMoore},
\cite{FeehanLeness}, \cite{PidstrigachTyurin}, and \cite{GottscheZagier}. 
\end{remark}

Theorem~\ref{thm:StrongAdjunctionInequality} follows from a
relation which holds for embedded surfaces with arbitrary
self-intersection number. This relation can be viewed as a
generalization of the relation appearing in~\cite{SympThom}. Once
again, we begin by stating the case when $b_1(X)=0$, in the interest
of exposition.

\begin{theorem}
\label{thm:Relation}
Let $X$ be a smooth, closed, connected, oriented four-manifold with $b_1(X)=0$,
and let
$\Sigma\subset X$ be a smoothly embedded surface with genus
$g(\Sigma)>0$. 
Then, for each
$\SpinC$ structure $\spinc$ with
$$-\langle c_1(\spinc),[\Sigma]\rangle+[\Sigma]\cm[\Sigma]\geq 0$$
and
$$
-\langle c_1(\spinc),[\Sigma]\rangle + [\Sigma]\cm [\Sigma] +
2d(\spinc) > 2g(\Sigma)-2,
$$
we have 
$$\SW_{X,\spinc}(U^{d}) =
\SW_{X,\spinc-\PD[\Sigma]}(U^{d'}),$$
where $d$ and $d'$ denote the dimensions of $\spinc$ and
$\spinc-\PD[\Sigma]$ respectively.
In the case where $b_2^+(X)=1$, both invariants are to be calculated
in the same component of $\Chambers(X)$.
\end{theorem}

More generally, we have the following.

\begin{theorem}
\label{thm:StrongRelation}
Let $X$ be a smooth, closed, connected, oriented four-manifold with $b_2^+(X)>0$.
Let $\Sigma\subset X$ be a surface with genus $g(\Sigma)>0$. Let $\ell$
be an integer so that there is a symplectic basis
$\{A_j,B_j\}_{j=1}^g$ for $H_1(\Sigma)$ so that $i_*(A_j)=0$ in
$H_1(X;\R)$ for $i=1,...,\ell$.  
For each $\SpinC$ structure $\spinc$ with
$$-\langle c_1(\spinc),[\Sigma]\rangle+[\Sigma]\cm[\Sigma]\geq 0$$
and each  $b\in\Alg(\Sigma)$ of degree $d(b)\leq \ell$ with
\begin{equation}
\label{eq:StrongAdjunctionInequalityFalse}
-\langle c_1(\spinc),[\Sigma]\rangle  + [\Sigma]\cm [\Sigma] +
2d(b) > 2g(\Sigma)-2,
\end{equation}
there is an element $b'\in\Alg(\Sigma)$ with $d(b')\geq d(b)$ so
that for any $a\in\Alg(X)$, we have
\begin{equation}
\label{eq:StrongRelation}
\SW_{X,\spinc}(a\cm i_*(b))=
\SW_{X,\spinc-\PD[\Sigma]}(a\cm i_*(b')).
\end{equation}
Furthermore, if $b=U^{d/2}$, then $b'-U^{d'/2}$ lies in the ideal
generated by $H_1(\Sigma)$ in $\Alg(\Sigma)$.  Once again, in the case
where $b_2^+(X)=1$, both invariants are to be calculated in the same
component of $\Chambers(X)$.
\end{theorem}

Theorem~\ref{thm:StrongAdjunctionInequality} is a simple consequence
of
Theorem~\ref{thm:StrongRelation}, as the following proof shows.

\vskip0.3cm
\noindent{\bf Theorem~\ref{thm:StrongRelation} $\Rightarrow$
Theorem~\ref{thm:StrongAdjunctionInequality}.}  Suppose
Theorem~\ref{thm:StrongAdjunctionInequality} were false; i.e. suppose
there were $X$, $\Sigma$, $\spinc$, $a$, and $b$ which satisfy the
hypotheses of the theorem, but which violate Adjunction
Inequality~\eqref{eq:StrongAdjunctionInequality}.  We can assume
without loss of generality that $$-\langle c_1(\spinc),[\Sigma]\rangle
+ [\Sigma]\cm[\Sigma] \geq 0,$$ by reversing the orientation of
$\Sigma$ if necessary (when $b_2^+(X)>1$). Thus,
Theorem~\ref{thm:StrongRelation} applies. Let $b'$ be the element
which satisfies Relation~\eqref{eq:StrongRelation}, so we have that
$\SW_{X,\spinc-\PD[\Sigma]}(a\cm i_*(b'))\neq 0$. Since $d(b)$ and
$d(b')$ are homogeneous elements with the same degree modulo two, and
$d(b')\geq d(b)$, it follows that we can find elements $a'\in\Alg(X)$
and $b''\in
\Alg(\Sigma)$ with $d(b'')=d(b)$, and
$\SW_{X,\spinc-\PD[\Sigma]}(a'\cm i_*(b''))\neq 0$. Now, since
$[\Sigma]\cm[\Sigma]\geq 0$ and $d(b'')=d(b)$, we see that $\Sigma$
also violates the adjunction inequality for $\spinc-\PD[\Sigma]$,
$a'\in\Alg(X)$, and $b''\in\Alg(\Sigma)$.  Proceeding in this way, we
see that $\spinc-n\PD[\Sigma]$ is a Seiberg-Witten basic class for all
$n\geq 0$. If $b_2^+(X)>1$, then there are only finitely many basic
classes of $X$, so since $\Sigma$ is not a torsion class, we get a
contradiction, proving Theorem~\ref{thm:StrongAdjunctionInequality} in
this case.

The above argument works also when $b_2^+(X)=1$, since there are still
only finitely many basic classes of the form $\spinc-n\PD[\Sigma]$ in
the chamber corresponding to $\PD[\Sigma]$.  We see this as
follows. Fix a metric $g$ on $X$ and a generic self-dual two-form
$\eta$. Clearly, if $\spinc$ is fixed and $n$ is sufficiently large,
the sign of $\PD[\Sigma]\cm\omega_g$ agrees with the sign of $-2\pi
c_1(\spinc-n\PD[\Sigma])\cm \omega_g +\int \eta\wedge \omega_g$;
i.e. for all large $n$, the $\eta$-perturbed moduli spaces for
$\spinc-n\PD[\Sigma]$ can be used calculate the invariant in the
component which contains $\PD[\Sigma]$. But the usual compactness
argument shows that all but finitely many of these moduli spaces are
empty.  Again, we have the contradiction completing the proof of
Theorem~\ref{thm:StrongAdjunctionInequality}. \hfill $\Box$
\vskip0.3cm

By blowing up, Theorem~\ref{thm:StrongRelation} is reduced to the case
where the self-intersection number of $\Sigma$ is sufficiently
negative.  The theorem is then proved by expressing the Seiberg-Witten
invariants of a four-manifold with such an embedded surface $\Sigma$
in terms of relative invariants, which take values in a
Seiberg-Witten-Floer homology associated to non-trivial circle bundles
over $\Sigma$. In the presence of the topological hypotheses on the
inclusion of $H_1(\Sigma)$ in $H_1(X)$, the above relation then
follows from properties of this Floer homology.

The outline of this paper is as follows. In
Section~\ref{sec:Examples}, we give examples which show that the
adjunction inequalities are sharp. Our examples include four-manifolds
with $b_2^+(X)=1$, and also examples where both $b_2^+(X)>1$ and
$b_1(X)>0$.  In Section~\ref{sec:RelProd}, we show how
Theorem~\ref{thm:StrongRelation} can be deduced from properties of a
product formula, which relates the Seiberg-Witten invariants of a
four-manifold containing an embedded surface with sufficiently
negative self-intersection number with certain relative invariants
associated to $X-\Sigma$.  For completeness, we also show how a
modified version of Theorem~\ref{thm:StrongRelation} implies
Theorem~\ref{thm:GeneralAdjunctionInequality}.  In
Section~\ref{sec:ModY} we review the gauge theory for circle bundles
over Riemann surfaces as developed in~\cite{MOY}.  There is one
$\SpinC$ structure in which the moduli space of reducibles has
singularities (to which we return in a later section). In
Section~\ref{sec:ProductFormula}, we prove the product formula
introduced in Section~\ref{sec:RelProd}, assuming technical facts
about the moduli spaces over $\NSig$, the tubular neighborhood of
$\Sigma$.  In Section~\ref{sec:SWirr}, we define an invariant with
irreducible boundary values and use properties of this relative
invariant to analyze the terms appearing in the product formula,
completing the proof of Theorem~\ref{thm:StrongRelation}. In
Section~\ref{sec:ModNSig}, we prove the technical facts about the
moduli spaces over $\NSig$ which were used in earlier sections. In
Section~\ref{sec:Perturbation}, we show how to extend the results of
Sections~\ref{sec:ModY} and
\ref{sec:ModNSig} to deal with the remaining $\SpinC$
structure. Finally, in Section~\ref{sec:Cohomology}, which should be
viewed as an appendix, we discuss representatives for the cohomology
classes used throughout the paper.

{\bf Acknowledgements.} The authors wish to thank Vicente Mu\~noz for
his very helpful comments on an early version of this paper.
\section{Examples}
\label{sec:Examples}

We give some examples now of four-manifolds $X$ which admit basic
classes of non-zero dimension. We begin by giving examples where
$b_2^+(X)>1$ and $b_1(X)>0$, to show that the adjunction
inequality in Theorem~\ref{thm:StrongAdjunctionInequality} is
sharp. (It is an open problem whether manifolds with $b_2^+(X)>1$ and
$b_1(X)=0$ can admit basic classes of non-zero dimension.)

\subsection{Examples of Theorem~\ref{thm:StrongAdjunctionInequality}
with $b_2^+(X)>1$}
\label{subsec:BPlusBig}

To construct these examples, we use the following
construction.

\begin{defn}
Let $X$ be smooth four-manifold and let $S\subset X$ be an embedded
two-sphere with zero self-intersection number.  Let $X'$ denote a
manifold obtained as surgery on $S$; i.e.  $$X' = (X -
\nbd{S})\cup_{\phi} S^1\times D^3,$$ where $\nbd{S}$ is an open tubular neighborhood of $S$ and
$\phi\colon \partial (X-\nbd{S})\longrightarrow S^1\times S^2$ is a
orientation-reversing diffeomorphism. Note that up to isotopy there
are two possible choices for $\phi$. Let $C\subset X'$ denote the
closed curve which is the core of the added $S^1\times D^3$. Note that
there is a diffeomorphism $X-S
\cong X'-C$.
\end{defn}

\begin{prop}
\label{prop:Surgery}
Let $X$ be a closed, smooth, oriented four-manifold with $b_2^+(X)>1$,
and let $S\subset X$ be a homologically trivial embedded
two-sphere.
For each $\SpinC$ structure $\spinc$ on $X$, there is a unique
induced $\SpinC$ structure $\spinc'$ on $X'$ with the property that
$$\spinc|_{X-S}=\spinc'|_{X'-C}.$$ Then, $d(\spinc')=d(\spinc)+1$; and
for all $a\in \Alg(X)$ $$\SW_{X',\spinc'}(a\cm \mu(C)) =
\SW_{X,\spinc}(a), $$
for some homology orientation on $X'$.
\end{prop}

\begin{proof}
The dimension statement is straightforward.

To prove the relation, we pull $X$ apart along $S^1\times
S^2=\partial~\nbd{S}$, and study the corresponding moduli spaces (see
Section~\ref{sec:ProductFormula} for more discussion on such matters).
Let $\Xtrunc$ denote the complement $X-S$, given a cylindrical-end
metric modeled on the product metric $[0,\infty)\times S^1\times S^2$,
where $S^2$ is given its standard, round metric. Note that this metric
can be extended over both $S^1\times D^3$ and $D^2\times S^2$ to give
metrics with non-negative scalar curvature. Consequently, the moduli
spaces of solutions over $S^1\times S^2$, $S^1\times D^3$, and
$D^2\times S^2$ consist entirely of smooth reducibles (i.e. the moduli
spaces are identified with $S^1$, $S^1$, and a point respectively).

Let $\ModSW_{\Xtrunc}(\spinc_0)$ denote the moduli space of finite
energy solutions to the Seiberg-Witten equations over $\Xtrunc$ in the
$\SpinC$ structure $\spinc_0=\spinc|_{\Xtrunc}$. Thus, we can think of the
boundary map as a map $$\Restrict\colon
\ModSW_{\Xtrunc}(\spinc_0)\longrightarrow S^1.$$ Gluing theory gives a
diffeomorphism for all sufficiently large $T>0$:
$$\ModSW_{X(T)}(\spinc)\cong\Restrict^{-1}(x_0),$$ where $X(T)$
denotes the metric on $X$ with neck-length $T$ and $x_0\in S^1$
corresponds to the unique reducible on $S^1\times S^2$ which extends
to $D^2\times S^2$. Consequently, 
\begin{equation}
\label{eq:InvariantOfOriginal}
\SW_{X,\spinc}(a)=\langle
\ModSW_{\Xtrunc}(\spinc_0), \mu(a)\cup \mu(C)\rangle,
\end{equation}
since $\mu(C)$ is represented by the holonomy class around $C$ (see
Proposition~\ref{prop:MuCircle}).

Similarly, gluing gives a diffeomorphism of
$$\ModSW_{\Xtrunc}(\spinc_0)\cong \ModSW_{X'(T)}(\spinc'),$$ and
consequently 
\begin{equation}
\label{eq:InvariantOfSurgery}
\SW_{X',\spinc'}(a\cm C)=\langle
\ModSW_{\Xtrunc}(\spinc_0), \mu(a\cm C)\rangle.
\end{equation}
Together, Equations~\eqref{eq:InvariantOfOriginal} and
\eqref{eq:InvariantOfSurgery} prove the proposition.
\end{proof}

\begin{remark}
Of course, the above result also holds when $b_2^+(X)=1$, provided
that both invariants are evaluated in the same chamber. 
\end{remark}

Now we construct our examples.
Fix natural numbers $n$, $k$, and $m$ with $2k\geq n>1$,
and let $X$ be the four-manifold $E(n)\# m (S^3\times S^1)$, where
$E(n)$ is a simply-connected elliptic surface with no multiple fibers
and with geometric genus $n-1$.  Let $\Sigma_0\subset E(n)$ denote a
symplectic submanifold representing the homology class $S + k F$,
where $S$ and $F$ denote the homology classes of a section and a fiber
respectively of the elliptic fibration.  Let $T_i\subset X$ denote a fiber in the
elliptic fibration of the $i^{th}$ summand $S^3\times S^1$.  Let
$\Sigma\subset X$ denote the internal connected sum of $\Sigma_0\# F_1
\#...\#F_m$. Note that $g(\Sigma)=k+m$ and $\Sigma\cm\Sigma=2k-n\geq 0$.
Let $\spinc$ be the $\SpinC$ structure over $X$ induced from the
canonical $\SpinC$ structure on $E(n)$, and let $b=\Alg(\Sigma)$ be
the product $B_1\cm...\cm B_m$ where $B_i\in H_1(X)$ generates $H_1$
of the $i^{th}$ copy of $S^1\times S^3$.  Note that $d(b)=m$ and
$\Sigma$ has a symplectic basis $\{A_i,B_i\}_{i=1}^{k+m}$ for which
$A_1,...,A_{m}$ are homologically trivial in $X$. By
Proposition~\ref{prop:Surgery}, $$\SW_{X,\spinc}(B_1\cm...\cm
B_m)=1,$$ so the data $X$, $b$, $\Sigma$, $\spinc$ satisfy the hypotheses of
Theorem~\ref{thm:StrongAdjunctionInequality}.  In fact, we see that
$$\Sigma\cm \Sigma + \langle c_1(\spinc),[\Sigma] \rangle + 2
d(b) = 2g(\Sigma)-2,$$ which shows that the inequality of the
theorem is sharp, for all choices of $g(\Sigma)> 0$,
$\Sigma\cm\Sigma\geq 0$, and $d(b)$.

\subsection{Ruled Surfaces: The Homological Hypotheses on $H_1(\Sigma)$}
\label{subsec:RSurfEx} 
By looking at ruled surfaces, we show that
Inequality~\eqref{eq:GeneralAdjunctionInequality} is sharp, and hence
that some homological hypotheses are necessary for the stronger
inequality (which appears in
Theorem~\ref{thm:StrongAdjunctionInequality}) to hold.

As mentioned before, one cannot hope for the adjunction inequality of
Theorem~\ref{thm:StrongAdjunctionInequality} to be valid without
additional topological hypotheses on the inclusion of $\Sigma$ in $X$.
Indeed, fix $n\geq 0$ and $g>0$, and let $X$ be the two-sphere bundle
over a surface $\Sigma$ of genus $g$, associated to the circle bundle
with Euler number $n$. In particular, $X$ contains an embedded copy of
$\Sigma$ with $\Sigma\cm\Sigma=n$. In the chamber corresponding to
$\PD[\Sigma]$, there is a zero-dimensional basic class 
$\spinc_0$ with $c_1(\spinc_0)=-\Canon{X}$, where $\Canon{X}$ is the
canonical class of $X$ viewed as K\"ahler manifold. Moreover, letting
$F$ be the class of the two-sphere fiber in $X$, we see that the
moduli space associated to $\spinc_0+d\PD[F]$ is identified with
$\Sym^d(\Sigma)$, and $U$ is the symmetric product of the volume form
of $\Sigma$ (see Proposition~\ref{prop:MuClassesOnSym} for a related
discussion). Thus, if $\spinc=\spinc_0+d\PD[F]$, then
$\SW_{\spinc}(U^d)\neq 0$, and $$\langle c_1(\spinc), [\Sigma]\rangle = 2d.$$
Clearly, Adjunction
Inequality~\ref{thm:GeneralAdjunctionInequality} is sharp for all
values of $k$, $d$, $n$, and $g$ provided that $-n\leq k$, where
$k=\langle c_1(\spinc),[\Sigma]\rangle$, $2d=d(b)$,
$n=\Sigma\cm\Sigma$, and $g=g(\Sigma)$.  (This construction, strictly
speaking, only gives us even values of $d(b)$. For odd values, one
can attach an $S^1\times S^3$.) In particular, we see that some
homological criterion on the embedding of $\Sigma\subset X$ is
necessary for the stronger
Inequality~\eqref{eq:StrongAdjunctionInequality} to hold.

\section{From Product Formulas to Relations}
\label{sec:RelProd}

The aim of this section is to outline the proof of
Theorem~\ref{thm:StrongRelation}. By employing the blowup formula in a
manner analogous to~\cite{SympThom}, we reduce to the case where the
self-intersection number of $\Sigma$ is very negative
(Proposition~\ref{prop:StrongRelationVeryNegative}). The invariants in
this latter case are studied via a product formula, which we state
(and prove in Section~\ref{sec:ProductFormula}), whose terms are then related  with other
Seiberg-Witten invariants of $X$. In the end of the section, we
discuss the modifications which are needed to prove
Theorem~\ref{thm:GeneralAdjunctionInequality}. 

We reduce Theorem~\ref{thm:StrongRelation} to the following special
case.

\begin{prop}
\label{prop:StrongRelationVeryNegative}
Theorem~\ref{thm:StrongRelation} holds, under the additional
hypothesis that
\begin{eqnarray*}
0\leq -\langle c_1(\spinc),[\Sigma]\rangle + [\Sigma]\cm[\Sigma]\leq 2g(\Sigma)-2
&{\text{and}}& -[\Sigma]\cm[\Sigma]> 2g-2.
\end{eqnarray*}
\end{prop}

The reduction involves the following basic result of Fintushel and
Stern.

\begin{theorem}(Blowup Formula) \cite{FS} and \cite{Salamon}.
Let $X$ be a smooth, closed four-manifold, and let ${\widehat
X}=X\#\CPbar$ denote its blow-up, with exceptional class $E\in
H^2({\widehat X};\Z)$. If $b_2^+(X)>1$, then for each $\SpinC$
structure ${\widehat\spinc}$ on ${\widehat X}$ with
$d({\widehat\spinc})\geq 0$, and each $a\in\Alg(X)\cong\Alg({\widehat
X})$, we have $$SW_{\widehat{X},{\widehat\spinc}}(a) =
SW_{X,\spinc}(U^{m} a),$$ where $\spinc$ is the $\SpinC$ structure
induced on $X$ obtained by restricting ${\widehat\spinc}$, and
$2m=d(\spinc)-d({\widehat
\spinc})$. If $b_2^+(X)=1$, there is a one-to-one correspondence
between components of $\Omega^+(X)$ and $\Omega^+({\widehat X})$, and
the above relation holds provided
both invariants are calculated in chambers associated to
corresponding components.
\end{theorem}

Before showing how to reduce Theorem~\ref{thm:StrongRelation} to the
special case, we point out that another special case of
Theorem~\ref{thm:StrongRelation} was already proved in Theorem
1.3~\cite{SympThom}. More specifically, the following was shown: \newline

\begin{theorem}~\cite{SympThom} 
\label{thm:SympThomRel}
Let $X$ be a smooth, closed, connected, oriented four-manifold with $b_2^+(X)>0$.
Let $\Sigma\subset X$ be a surface with genus $g(\Sigma)>0$ and
negative self-intersection. 
For each $\SpinC$ structure $\spinc$ with
$$
-\langle c_1(\spinc),[\Sigma]\rangle  + [\Sigma]\cm [\Sigma]  > 2g(\Sigma)-2,
$$
there is an element $b'\in\Alg(\Sigma)$ so that
that for any $a\in\Alg(X)$, we have
$$
\SW_{X,\spinc}(a)
=
\SW_{X,\spinc-\PD[\Sigma]}(a\cm i_*(b')).
$$
Furthermore, $b'-U^{d'/2}$ lies in the ideal generated by
$H_1(\Sigma)$ in $\Alg(\Sigma)$. 
\end{theorem}

\begin{remark}
In the language of Theorem~\ref{thm:StrongRelation}, this case
corresponds to $\ell=0$ and $b=1$. 
\end{remark}

\vskip0.3cm
\noindent{\bf Proposition~\ref{prop:StrongRelationVeryNegative}
$\Rightarrow$ Theorem~\ref{thm:StrongRelation}.}
Let $g=g(\Sigma)$, 
fix an integer $m$ with 
$$m > [\Sigma]\cm[\Sigma] + 2g-2,$$
let ${\widehat X}=X\#m\CPbar$, and let ${\widehat \Sigma}$ be the
``proper transform'' of $\Sigma$, the embedded surface obtained by
internal connected sum of $\Sigma$ with the $m$ exceptional spheres in
the $\CPbar$ summands; i.e.
$$\PD[{\widehat \Sigma}]=\PD[\Sigma]-E_1-...-E_m.$$
Finally, let ${\widehat\spinc}$ denote the $\SpinC$ structure on
${\widehat X}$ which agrees with $\spinc$ in the complement of the
exceptional spheres, whose Chern class satisfies
$$c_1({\widehat\spinc})=c_1(\spinc)-E_1-...-E_m.$$
It is easy to check that:
\begin{eqnarray*}
-[\widehat\Sigma]\cm[\widehat\Sigma]=m-[\Sigma]\cm[\Sigma]&>&
2g-2,\\
-\langle c_1({\widehat\spinc}),[{\widehat\Sigma}]\rangle
+[{\widehat \Sigma}]\cm[{\widehat \Sigma}]
&=& -\langle c_1(\spinc),[\Sigma]\rangle+[\Sigma]\cm[\Sigma];
\\
d(\spinc) &=& d({\widehat\spinc}).
\end{eqnarray*}
Now, if $$-\langle c_1(\spinc),[\Sigma]\rangle+[\Sigma]\cm[\Sigma]\geq
2g,$$ the hypotheses of Theorem~\ref{thm:SympThomRel} are satisfied;
and otherwise, the hypotheses of
Proposition~\ref{prop:StrongRelationVeryNegative} are. In either case,
for each $b\in\Alg(\Sigma)$ of degree $d(b)\leq \ell$,
we can find $b'\in\Alg(\Sigma)$ with
\begin{equation}
\label{eq:RelationOnBlowup}
\SW_{{\widehat X},{\widehat \spinc}}(a\cm i_*(b))=
\SW_{{\widehat X},{\widehat
\spinc}-\PD[{\widehat \Sigma}]}(a\cm i_*(b')).
\end{equation}
According to the  blow-up formula, 
\begin{equation}
\label{eq:BlowupOne}
\SW_{{\widehat X},{\widehat \spinc}}(a\cm i_*(b))=
\SW_{X,\spinc}(a\cm i_*(b));
\end{equation}
and,
since
${\widehat \spinc}-\PD[{\widehat \Sigma}]$ agrees with $\spinc-\PD[\Sigma]$
away from the exceptional spheres and
$$c_1({\widehat \spinc}-\PD[{\widehat
\Sigma}])=c_1(\spinc-\PD[\Sigma])-E_1-...-E_m,$$
we see from another application of the blowup formula that
\begin{equation}
\label{eq:BlowupTwo}
\SW_{{\widehat X},{\widehat
\spinc}-\PD[{\widehat \Sigma}]}(a\cm i_*(b'))=
\SW_{X,\spinc-\PD[\Sigma]}(a\cm i_*(b')).
\end{equation}
Theorem~\ref{thm:StrongRelation} then follows by combining
Equations~\eqref{eq:RelationOnBlowup}, \eqref{eq:BlowupOne} and
\eqref{eq:BlowupTwo}.\hfill $\Box$
\vskip0.3cm

We now turn to the special case considered in
Proposition~\ref{prop:StrongRelationVeryNegative}.  We will study the
Seiberg-Witten invariant of $X$ by decomposing it into two pieces
$$X=\NSig\cup_Y(X-\NSig),$$ where $Y$ a circle bundle over $\Sigma$
(as in the proposition), and $\NSig$ is the associated disk bundle.
Following~\cite{MOY}, the moduli space of Seiberg-Witten monopoles
over $Y$ decomposes into an irreducible and a reducible
component. (Actually, there is one $\SpinC$ structure over $Y$, where
it is necessary to perturb the equations for this decomposition to
occur; this perturbation is studied Section~\ref{sec:Perturbation}.)
Correspondingly, we construct relative invariants of $X-\Sigma$,
denoted $\SWirr_\spinc$ and $\SWred_\spinc$, arising from the $L^2$
moduli spaces on $X-\Sigma$ with irreducible and reducible boundary
values.  In Section~\ref{sec:ProductFormula} (see
Lemma~\ref{lemma:TwoCases}, and the discussion following it), we
prove the following:

\begin{prop}
\label{prop:ProductFormula}
Suppose
\begin{eqnarray*}
0\leq -\langle c_1(\spinc),[\Sigma]\rangle + [\Sigma]\cm[\Sigma]\leq 2g(\Sigma)-2
&{\text{and}}& -[\Sigma]\cm[\Sigma]> 2g-2.
\end{eqnarray*}
Then,
$$\SW_{X,\spinc}=\SWirr_{\spinc} + \SWred_{\spinc}.$$
\end{prop}

We can interpret the latter invariant in terms of the closed manifold
as follows.

\begin{defn}
\label{def:xis}
Let $\Sigma$ be a surface of genus $g$, and let $\{A_i,B_i\}_{i=1}^g$
be a standard symplectic basis for $H_1(\Sigma;\Z)$.  For $j=0,...,g$,
let $\xi_j([\Sigma])\in\Alg(\Sigma)$ be the degree $2j$ component of
$$\prod_{i=1}^g\Big(1+U+A_i\cm B_i\Big)\in\Alg(\Sigma);$$
i.e. $\xi_0=1$, $\xi_1(\Sigma)=gU+\sum A_i\cm B_i$, ...,
$\xi_g(\Sigma)=\prod_{i=1}^g\Big(U+A_i\cm B_i\Big)$.
\end{defn}

\begin{prop}
\label{prop:RedBoundary}
Suppose
\begin{eqnarray*}
0\leq -\langle c_1(\spinc),[\Sigma]\rangle + [\Sigma]\cm[\Sigma]\leq 2g(\Sigma)-2
&{\text{and}}& -[\Sigma]\cm[\Sigma]> 2g-2.
\end{eqnarray*}
Then, letting $$e=g-1+\frac{\langle
c_1(\spinc),[\Sigma]\rangle-[\Sigma]\cm[\Sigma]}{2},$$ we have that
$$\SWred_{\spinc}(a)=\SW_{X,\spinc-\PD[\Sigma]}(a\cm
\xi_{g-1-e}([\Sigma]))$$ for all $a\in\Alg(X)$.
\end{prop}

Furthermore, under the homological condition of
Theorem~\ref{thm:StrongRelation}, we will express $\SWirr_{\spinc}$ in
terms of $\SW_{X,\spinc-\PD[\Sigma]}$, as follows.

\begin{prop}
\label{prop:Vanishing}
Suppose 
\begin{eqnarray*}
0\leq -\langle c_1(\spinc),[\Sigma]\rangle + [\Sigma]\cm[\Sigma]\leq 2g(\Sigma)-2&{\text{and}}&
-[\Sigma]\cm[\Sigma]> 2g-2,
\end{eqnarray*}
and let $\ell$ be an integer so that there is a symplectic basis
$\{A_i,B_i\}_{i=1}^g$ for $H_1(\Sigma)$ so that $i_*(A_i)=0$ in
$H_1(X;\R)$ for $i=1,...,\ell$.  Then, for each $b\in\Alg(\Sigma)$ of
degree $e<d(b)\leq\ell$, there is an element $b''\in\Alg(\Sigma)$
so that
$$\SWirr_{\spinc}(a\cm i_*(b))=
\SW_{\spinc-\PD[\Sigma]}(a\cm i_*(b'')).$$
Furthermore, $b_2$ lies in the ideal generated by
$H_1(\Sigma)$ in $\Alg(\Sigma)$.
\end{prop}

The proof of Proposition~\ref{prop:Vanishing} is given in the end of
Section~\ref{sec:SWirr}. 

Proposition~\ref{prop:StrongRelationVeryNegative} follows immediately
from Propositions~\ref{prop:ProductFormula}--\ref{prop:Vanishing}.  In
the proof of these latter propositions, we will construct a natural
Seiberg-Witten-Floer functor for four-manifolds which bound $Y$.

Before proceeding, we pause to tie up one more loose end:
Theorem~\ref{thm:GeneralAdjunctionInequality}.  That result can be
reduced to a relation which replaces Theorem~\ref{thm:StrongRelation},
using the same argument given in the proof of
Theorem~\ref{thm:StrongAdjunctionInequality}.  The relevant relation
in this case is:

\begin{theorem}
\label{thm:GeneralRelation}
Let $X$ be a smooth, closed, connected, oriented four-manifold with $b_2^+(X)>0$.
Let $\Sigma\subset X$ be a surface with genus $g(\Sigma)>0$. 
For each $\SpinC$ structure $\spinc$ with
$$-\langle c_1(\spinc),[\Sigma]\rangle+[\Sigma]\cm[\Sigma]\geq 0$$
and each  $b\in\Alg(\Sigma)$ of degree $d(b)$ with
\begin{equation}
\label{eq:GeneralAdjunctionInequalityFalse}
-\langle c_1(\spinc),[\Sigma]\rangle  + [\Sigma]\cm [\Sigma] +
d(b) > 2g(\Sigma)-2,
\end{equation}
there is an element $b'\in\Alg(\Sigma)$ with $d(b')\geq d(b)$ so
that for any $a\in\Alg(X)$, we have
\begin{equation}
\label{eq:GeneralRelation}
\SW_{X,\spinc}(a\cm i_*(b))
=
\SW_{X,\spinc-\PD[\Sigma]}(a\cm i_*(b')).
\end{equation}
Furthermore, if 
$b=U^{d}$, then $b'-U^{d'}$ lies in the ideal generated by
$H_1(\Sigma)$ in $\Alg(\Sigma)$. 
\end{theorem}

Once again, via the blowup formula,
this relation can be reduced to the case where the self-intersection
number $\Sigma$ is very negative; i.e. 
\begin{eqnarray*}
0\leq -\langle c_1(\spinc),[\Sigma]\rangle + [\Sigma]\cm[\Sigma]\leq 2g(\Sigma)-2
&{\text{and}}& -[\Sigma]\cm[\Sigma]> 2g-2.
\end{eqnarray*}
(compare Proposition~\ref{prop:StrongRelationVeryNegative}).
Like Proposition~\ref{prop:StrongRelationVeryNegative}, this
special case also follows from the product formula in
Proposition~\ref{prop:ProductFormula}, the relation in
Proposition~\ref{prop:RedBoundary}, together with
the following analogue of Proposition~\ref{prop:Vanishing} (whose
proof is also given in the end of Section~\ref{sec:SWirr}):

\begin{prop}
\label{prop:GeneralVanishing}
Suppose 
\begin{eqnarray*}
0\leq -\langle c_1(\spinc),[\Sigma]\rangle + [\Sigma]\cm[\Sigma]\leq 2g(\Sigma)-2&{\text{and}}&
-[\Sigma]\cm[\Sigma]> 2g-2.
\end{eqnarray*}
Then, for each $b\in\Alg(\Sigma)$ of degree $2e<d(b)$, there is an
$b''\in\Alg(\Sigma)$ so that 
$$
\SWirr_{\spinc}(a\cm i_*(b))
=
\SW_{\spinc-\PD[\Sigma]}(a\cm i_*(b'')).
$$ Furthermore, $b''$ lies
in the ideal generated by $H_1(\Sigma)$.
\end{prop}

\section{Gauge theory on  $\R\times Y$}
\label{sec:ModY}

The Seiberg-Witten moduli spaces over $Y$ and $\R\times Y$ were
studied for Seifert fibered three-manifolds $Y$ in~\cite{MOY}. We
summarize these results here, for $Y$ a circle-bundle over a Riemann
surface $\Sigma$ with $g(\Sigma)>0$ and Euler number $-n$, where
$n> 2g-2$.

$Y$ admits a canonical $\SpinC$ structure
whose bundle of spinors is
$\C\oplus\pi^*(\Canon{\Sigma}^{-1})$, which we use to
identify the $\SpinC$ structures on $Y$ with $H^2(Y;\Z)\cong \Z^{2g}\oplus
\Z/n\Z$. 

Let $\ModSWThree_Y(\spinct)$ denote the moduli space of solutions to
the Seiberg-Witten equations over $Y$ in the $\SpinC$ structure
$\spinct$. Here, we use the metric $g_Y$ and $SO(3)$-connection over $TY$
of~\cite{MOY}.  Given a pair of components $C_1,C_2$ in
$\ModSWThree_Y(\spinct)$, let $\ModFlow(C_1,C_2)$ denote the moduli
space of solutions $[A,\Phi]$ to the Seiberg-Witten equations on
$\R\times Y$ for which
\begin{eqnarray*}
\lim_{t\goesto -\infty}[A,\Phi]|_{\{t\}\times Y}\in C_1, &{\text{and}}&
\lim_{t\goesto \infty}[A,\Phi]|_{\{t\}\times Y}\in C_2.
\end{eqnarray*}
This moduli space admits a translation action by $\R$. Let
$\UnparamModFlow(C_1,C_2)$ denote the quotient of $\ModFlow(C_1,C_2)$ by the
translation action.

In general, these spaces admit a Morse-theoretic interpretation.  If
$c_1(\spinct)$ is a torsion class, there is a real-valued functional
$$\CSD\colon
\SWConfig(Y,\spinct)\longrightarrow
\R$$ 
defined over the configuration space $\SWConfig(Y,\spinct)$ of pairs
$(B,\Psi)$ of spin-connections $B$ in $\spinct$ and spinors $\Psi$
modulo gauge.  The critical manifolds are the moduli spaces
$\ModSWThree(Y;\spinct)$.  When $c_1(\spinct)$ is not torsion, the
functional is circle-valued.  The Seiberg-Witten equations on
$\R\times Y$ are the upward gradient-flow equations for this
functional. In keeping with this interpretation, we call
$\UnparModFlow(C_1,C_2)$ the space of unparameterized flows
from $C_1$ to $C_2$.

\begin{theorem}
\label{thm:MOY}
(\cite{MOY})
Let $Y$ be a circle-bundle over a Riemann surface with genus $g>0$ and
Euler number $-n< 2-2g$.  The moduli space $\ModSWThree_Y(\spinct)$ is
empty unless $\spinct$ corresponds to a torsion class in
$H^2(Y;\Z)$. So, suppose $\spinct$ corresponds to $\mult\in\Zmod{n}\subset
H^2(Y;\Z)$.
\begin{list}
{(\arabic{bean})}{\usecounter{bean}\setlength{\rightmargin}{\leftmargin}}
\item
\label{item:SmallK}
If $0\leq \mult <g-1$ then $\ModSWThree_Y(\spinct)$ contains two
components, a reducible one $\Jac$, identified with the Jacobian torus
$H^1(\Sigma;\R/\Z)$, and a smooth irreducible
component $\CritMan$ diffeomorphic to $\Sym^\mult(\Sigma)$. Both of these
components are non-degenerate in the sense of Morse-Bott. There is an
inequality $\CSD(\Jac)>\CSD(\CritMan)$, so
the space
$\UnparModFlow(\Jac,\CritMan)$ is empty.  The space
$\UnparModFlow(\CritMan,\Jac)$ is smooth of expected dimension $2\mult$;
indeed it is diffeomorphic to $\Sym^\mult(\Sigma)$.
\newline
\item If $g-1<\mult\leq 2g-2$, the Seiberg-Witten moduli spaces over both
$Y$ and $\R\times Y$
in this $\SpinC$ structure are naturally identified with the	
corresponding moduli spaces in the $\SpinC$ structure $2g-2-\mult$, which
we just described. 
\newline
\item For all other $\mult\neq g-1$, $\ModSWThree_Y(\spinct)$ contains only
reducibles. Furthermore, it is smoothly identified with the Jacobian
torus.
\end{list}
\end{theorem}

In the $\SpinC$ structure corresponding to $g-1\in\Zmod{n}$, the
unperturbed Seiberg-Witten equations used in Theorem~\ref{thm:MOY} are
inconvenient, since the corresponding reducible manifold is not smooth
in the sense of Morse-Bott.  To overcome this difficulty, when working
in this $\SpinC$ structure, we use a perturbation of the equations
where the theory resembles the case where $0\leq e<g-1$ (and, in
particular, the reducibles are smooth).  A thorough discussion
of the perturbation is given in Section~\ref{sec:Perturbation}.

\section{The Product Formula}
\label{sec:ProductFormula}


In this section, we define two quantities, $\SWirr$ and $\SWred$, and
prove that the Seiberg-Witten invariant decomposes into a sum of these
(Propostion~\ref{prop:ProductFormula}).  
Furthermore, we express $\SWred$ in
terms of another Seiberg-Witten invariant of $X$
(Proposition~\ref{prop:RedBoundary}).

Decompose $X$ as $$X=\NSig\cup_Y\Xtrunc,$$ where $Y$ is unit circle
bundle over $\Sigma$ with Euler number $-n$, with $n> 2g-2$.  $N$
is a tubular neighborhood of the surface $\Sigma$ (which is
diffeomorphic to the disk bundle associated to $Y$), and $\Xtrunc$ is
the complement in $X$ of the interior of $N$.  Fix metrics
$g_{\Xtrunc}$, $g_{N}$, and $g_Y$ for which $g_{\Xtrunc}$ and $g_{N}$
are isometric to $$dt^2 + g_Y^2$$ in a collar neighborhood of their
boundaries (where $t$ is a normal coordinate to the boundary). Let
$X(T)$ denote the Riemannian manifold which is diffeomorphic to $X$
and whose metric $g_T$ is obtained from the description
$$X(T)=\NSig\cup_{\partial\NSig=\{-T\}\times Y} [-T,T]\times
Y\cup_{\{T\}\times Y=-\partial \Xtrunc} \Xtrunc;$$
i.e. $g_T|_{N}=g_N$, $g_T|_{[-T,T]\times Y} = dt^2 + g_Y^2$, and
$g_T|_{\Xtrunc}=g_{\Xtrunc}$.  Our goal here is to provide, for all
sufficiently large $T$, a description of the moduli space
$\ModSW_{X(T)}(\spinc)$ on $X(T)$ in terms of the moduli spaces for
$Y$, $\ModSWThree_Y(\spinc|_Y)$, and the finite-energy,
cylindrical-end moduli spaces associated to $\Xtrunc$ and $N$, denoted
$\ModSW_{\Xtrunc}(\spinc|_{\Xtrunc})$, and $\ModSW_{\NSig}(\spinc|_N)$
respectively. In this context, finite energy means that the total
variation of the Chern-Simons-Dirac functional over the infinite
cylinder is bounded. Henceforth, $\Xtrunc$ and $N$ will denote the cylindrical-end
manifolds obtained by attaching $[0,\infty)\times Y$ (with appropriate
orientations) to the corresponding subsets of $X$.

In the case where $b_2^+(X)=1$, we choose the perturbing form $\eta$
to be compactly supported in $\Xtrunc$ in such a way that
$$-2\pi c_1(\spinc)\cm\omega_\infty+\int_{\Xtrunc}\eta\wedge
\omega_\infty$$ has the same sign as $\gamma\cm\omega_\infty$, where
$\gamma$ is a compactly supported representative for a class in the
chosen component $\Chamber\subset \Chambers(X)$, and $\omega_\infty$
is a self-dual harmonic two-form over $\Xtrunc$ with $\int_{\Xtrunc}
\omega_\infty\wedge \omega_\infty = 1$.  Note that such a $\gamma$ and
$\omega_\infty$ can be found since $\Sigma\cm\Sigma<0$, forcing
$b_2^+(\Xtrunc)=1$ (see~\cite{APSI}).  Now, the moduli spaces of the
$\eta$-perturbed Seiberg-Witten equations over $X(T)$ calculate the
invariant in the chosen chamber for all sufficiently large $T$.

We collect useful facts about the moduli spaces $\ModSW_N(\spinc|_N)$,
most of which we defer to Section~\ref{sec:ModNSig} (see
also~\cite{SympThom}), but first we introduce some notation. The map
$$\SpinC(N)\rightarrow \Z$$ given by $$\spinc\mapsto \langle
c_1(\spinc),[\Sigma]\rangle$$ induces a one-to-one correspondence
between $\SpinC$ structures and integers which are congruent to $n$
modulo $2$. Note that the $\SpinC$ structure over $Y$ $\spinc|_Y$
corresponds to the mod $n$ reduction of 
$$e=g-1+\frac{\langle c_1(\spinc),[\Sigma]\rangle +n}{2}$$
appearing in Theorem~\ref{thm:MOY}.

By taking limits at the end of the tube, one can define
maps
\begin{eqnarray*}
\Restrict\colon \ModSW_{\NSig}(\spinc)\longrightarrow
\ModSWthree_Y(\spinc|_Y)
&{\text{and}}&
\Restrict\colon \ModSW_{\Xtrunc}(\spinc)\longrightarrow
\ModSWthree_Y(\spinc|_Y)
\end{eqnarray*}
(see~\cite{MMR}).  If $\CritMan$ is a connected manifold of
$\ModSWThree_Y(\spinc|_Y)$, then $\ModSW_{\NSig}(\spinc,\CritMan)$ and
$\ModSW_{\Xtrunc}(\spinc,\CritMan)$ denotes the pre-image of
$\CritMan$ under $\Restrict$. Throughout the following discussion, we
will use the perturbation discussed in Section~\ref{sec:Perturbation}
over $\Xtrunc$, $\NSig$, and $Y$, when $\spinc|_Y$ corresponds to
$e=g-1$ (in the notation of Section~\ref{sec:ModY}); i.e. in this
case, $\ModSWThree_Y(\spinc|_Y)$, $\ModSW_{\NSig}(\spinc)$ and
$\ModSW_{\Xtrunc}(\spinc|_{\Xtrunc})$ will denote the perturbed
versions of these moduli spaces, with perturbation parameter $u$ in
the range $0<u<2$, in the notation of Section~\ref{sec:Perturbation}.
(We will show in Section~\ref{sec:Perturbation} that this is an allowable
perturbation to use when $b_2^+(X)=1$; i.e. we are computing the
Seiberg-Witten invariants in the correct chamber.)  When they are
clear from the context, we leave the $\SpinC$ structures out of the
notation.  Note that on the cylinders, the analogous boundary value
maps factor through the unparameterized spaces, defining
\begin{eqnarray*}
\Restrict_{\Jac}\colon \UnparModFlow(\CritMan,\Jac)\longrightarrow \Jac
&{\text{and}}&
\Restrict_{\CritMan}\colon
\UnparModFlow(\CritMan,\Jac)\longrightarrow \CritMan, 
\end{eqnarray*}
where $\Jac$ and $\CritMan$ are the critical manifolds of
Theorem~\ref{thm:MOY}. 

\begin{prop}
\label{prop:DimensionCounting}
Suppose that $\BigSpinc$,
and let $$e=g-1+\frac{\langle c_1(\spinc),[\Sigma]\rangle +n}{2},$$
Then according to Theorem~\ref{thm:MOY}, and
Theorem~\ref{thm:MOYPert} when $e=g-1$,
$\ModSWThree_Y(\spinc|_Y)$ has
two components, $\Jac$ and $\CritMan$, where $\CritMan$ is
diffeomorphic to $\Sym^\mult(\Sigma)$.
Furthermore, the expected dimensions of the moduli spaces over $N$ and $\Xtrunc$
are given by:
\begin{eqnarray}
\edim \ModSW_N(\Jac)&=& 2\mult+1 \\
\edim \ModSW_N(\CritMan)&=& 2\mult \\
\edim \ModSW_{\Xtrunc}(\Jac) &=& 2d+2g-2\mult-2 \\
\edim \ModSW_{\Xtrunc}(\CritMan) &=& 2d,
\end{eqnarray}
where $d=d(\spinc)$ and $g=g(\Sigma)$.
Moreover, $\ModIrr_N(\Jac)$, $\ModSW_N(\CritMan)$,
$\ModSW_{\Xtrunc}(\Jac)$, and $\ModSW_{\Xtrunc}(\CritMan)$ are
transversally cut out by the Seiberg-Witten equations (in particular,
they are manifolds of the expected dimension).
\end{prop}

\begin{proof}
This is a combination of Proposition~\ref{prop:NSigDefTheory} and 
\ref{prop:XtruncDefTheory} when $\langle c_1(\spinc),[\Sigma]\rangle\neq n$, 
and Proposition~\ref{prop:DimensionCountingPert} in the remaining
case.
\end{proof}

When studying the deformation theory of reducibles inside
$\ModSW_N(\Jac)$, the kernel and the cokernel of the Dirac operator
play a central role. These spaces can be concretely understood, thanks
to the holomorphic interpretation of the Dirac operator (see
also~\cite{SympThom}).

\begin{prop}
\label{prop:KerDirac}
Suppose that $\BigSpinc$,
then there is a natural correspondence between reducibles
$[(A,0)]\in\ModSW_N(\Jac)$ with holomorphic line bundles ${\mathcal E}$ 
of degree $\mult$ over $\Sigma$ which identifies
\begin{eqnarray*}
\Ker \Dirac_A = H^0(\Sigma,{\mathcal E}) &{\text{and}}& 	
\CoKer \Dirac_A = H^1(\Sigma,{\mathcal E}).
\end{eqnarray*}
\end{prop}

\begin{proof}
This follows from Theorem~\ref{thm:ExtendRed} and
Proposition~\ref{prop:RSurfCohomology} (see also the proof of
Theorem~\ref{thm:MOYPert} in the perturbed case).
\end{proof}

The above proposition allows us to understand an important class of
reducibles. 

\begin{defn}
The {\em jumping locus} $\Theta\subset\ModSW_N(\Jac)$ is the locus of
reducible solutions $[(A,0)]\in\ModSW_N(\Jac)$ for which $\Ker\Dirac_A$
is non-trivial.
\end{defn}

\begin{cor}
\label{cor:DimNSLocus}
Suppose that $\BigSpinc$, then
the jumping locus $\Theta\subset\Jac= \ModRed_N(\Jac)$ 
is the image of a smooth map
$\Sym^\mult(\Sigma)\longrightarrow \Jac$.
\end{cor}

\begin{proof}
According to Proposition~\ref{prop:KerDirac}, the space 
$\Theta\subset\Jac$ is identified with the space of degree
$\mult$ line bundles over $\Sigma$ with non-trivial $H^0$. The
forgetful map $\Sym^\mult(\Sigma)\longrightarrow \Jac$ which takes a
degree $\mult$ divisor, thought of as a complex line bundle with
section, to the underlying complex line bundle gives the surjection to
this locus.
\end{proof}

We will also need to understand those $\SpinC$ structures
$\spinc\in\SpinC(\NSig)$ for which $-n<
\langle c_1(\spinc),[\Sigma]\rangle\leq n$.

\begin{prop}
\label{prop:UnobstructedOnN}
If $$\SmallSpinc,$$ then the moduli
space $\ModSW_N(\Jac)$ contains only reducibles. Moreover, the space
of reducibles is smoothly identified with the Jacobian torus $\Jac$
(i.e. the kernel and the cokernel of the Dirac operator coupled to any
reducible vanishes). Furthermore, $\ModSW_N(\CritMan)$ is empty.
\end{prop}

\begin{proof}
When $|\langle c_1(\spinc),[\Sigma]\rangle|<  n$,
this is proved in Section~\ref{sec:ModNSig}, where it appears as
Proposition~\ref{prop:UnobstructedOnNNoPert}. The remaining case is
covered by Proposition~\ref{prop:UnobstructedOnNPert}.
\end{proof}

With these preliminaries in place, we turn to the Seiberg-Witten
invariants of $X$, by investigating the moduli spaces over $X(T)$.
Specifically, choose some $a\in\Alg(X)$ of degree $d(\spinc)$, and
indeed choose representatives for the corresponding homology classes
which are compactly supported in $\Xtrunc$. Let $\Divisors(a)$ denote
the corresponding representatives for $\mu(a)$ in the configuration
spaces for $\Xtrunc$ and $X(T)$ as appropriate (see
Section~\ref{sec:Cohomology} for a discussion of such
representatives). Recall that $\SW_{X,\spinc}(a)$ is the number of
points in $\ModSW_{X(T)}(\spinc)\cap \Divisors(a)$, counted with
appropriate sign.

\begin{lemma}
\label{lemma:TwoCases}
Suppose that  $\BigSpinc$, then
for each $\epsilon>0$, there is a $T_0>0$ so that for all $T\geq
T_0$ the restriction of $[(A,\Phi)]\in\ModSW_{X(T)}(\spinc)\cap \Divisors(a)$ to
any slice $\{t\}\times Y$ with $t\in[-T_0,T_0]$ lies within
$\epsilon$ (in the $\Cinfty$ topology) from either $\Jac$ or
$\CritMan$.  Accordingly, if $\epsilon$ is sufficiently small,
then $[(A,\Phi)]$ satisfies exactly one of the
following two conditions:
\begin{list}
{(H-\arabic{bean})}{\usecounter{bean}\setlength{\rightmargin}{\leftmargin}}
\item
\label{item:Red}
$[(A,\Phi)]|_{N}$ is $\Cinfty$ close to smooth reducible
and
$[(A,\Phi)]|_{\Xtrunc}$ is $\Cinfty$ close to (the restriction to
$\Xtrunc$)
of a configuration in  $\ModSW_{\Xtrunc}(\Jac)\cap\Divisors(a)$;
\newline
\item 
\label{item:Irr}
$[(A,\Phi)]|_{N}$ is $\Cinfty$ close to a configuration in
$\ModSW_{\NSig}(\CritMan)$, and
$[(A,\Phi)]|_{\Xtrunc}$ is $\Cinfty$ close to a configuration in
the cut-down moduli space $\ModSW_{\Xtrunc}(\CritMan)\cap \Divisors(a)$
\end{list}
\end{lemma}

\begin{proof}
This is a dimension-counting argument. Suppose we have a sequence
$[A_i,\Phi_i]\in \ModSW_{X(T_i)}(\spinc)\cap \Divisors(a)$, for some
increasing, unbounded sequence $\{T_i\}_{i=1}^{\infty}$ of real
numbers. By local compactness, there is a subsequence which converges
in $\CinftyLoc$ to a pair of configurations
$(A_N,\Phi_N)$ and $(A_{\Xtrunc},\Phi_{\Xtrunc})$ over $N$ and $\Xtrunc$
respectively. By the usual compactness arguments
(see~\cite{KMthom}), the total variation of the Chern-Simons-Dirac
functional of $(A_i,\Phi_i)$ over the cylinder $[-T_i,T_i]\times Y$
remains globally bounded (independent of $i$), so
$(A_N,\Phi_N)$ and $(A_{\Xtrunc},\Phi_{\Xtrunc})$ both have finite
energy.

First, we prove that either Hypothesis~(H-\ref{item:Red}) or
(H-\ref{item:Irr}) is satisfied.  There are {\em a priori} four cases,
according to which critical manifolds
$\Restrict[A_{\Xtrunc},\Phi_{\Xtrunc}]$ and $\Restrict[A_N,\Phi_N]$
lie in.
\begin{list}
{(P-\arabic{bean})}{\usecounter{bean}\setlength{\rightmargin}{\leftmargin}}
\item
\label{item:IrrRed}
The case where $\Restrict[A_N,\Phi_N]\in \Jac$ while
$\Restrict(A_{\Xtrunc},\Phi_{\Xtrunc})\in \CritMan$ is excluded
because $\CSD(\CritMan)>\CSD(\Jac)$.
\item
\label{item:RedIrr}
The case where $\Restrict[A_N,\Phi_N]\in \CritMan$ while
$\Restrict(A_{\Xtrunc},\Phi_{\Xtrunc})\in \Jac$ is excluded by a
dimension count, as follows. In this case, we see that
$\Restrict[A_{\Xtrunc},\Phi_{\Xtrunc}]\in
\Restrict_{\Jac}(\ModFlow(\CritMan,\Jac))\cap
\Restrict(\ModSW_{\Xtrunc}(\Jac)\cap \Divisors(a))$. But 
$$\Restrict_{\Jac}(\ModFlow(\CritMan,\Jac))=\Restrict_{\Jac}(\UnparModFlow(\CritMan,\Jac)),$$
so
\begin{eqnarray*}
\edim(\Restrict_{\Jac}(\ModFlow(\CritMan,\Jac))\cap
\Restrict(\ModSW_{\Xtrunc}(\Jac)\cap \Divisors(a)) &=& -2.
\end{eqnarray*}
It follows from Theorems~\ref{thm:MOY} and \ref{thm:MOYPert}
that $\ModFlow(\CritMan,\Jac)$ is smooth of the expected dimension, so
from the usual transversality results, the above intersection is
generically empty.
\item
\label{item:RedRed} 
Suppose that $\Restrict (A_N,\Phi_N) \in\Jac$ and $\Restrict
(A_{\Xtrunc},\Phi_{\Xtrunc})\in \Jac$. Then we see that
$$\Restrict[A_N,\Phi_N]=\Restrict[A_{\Xtrunc},\Phi_{\Xtrunc}]\in
\Restrict(\ModSW_N(\Jac))\cap
\Restrict(\ModSW_{\Xtrunc}(\Jac)\cap\Divisors(a));$$
but, according to Proposition~\ref{prop:DimensionCounting}
\begin{eqnarray*}
\lefteqn{\edim 
\Restrict(\ModIrr_N(\Jac))\cap
\Restrict(\ModSW_{\Xtrunc}(\Jac)\cap\Divisors(a))} \\
&=&  
\edim\ModSW_N(\Jac) + \ModSW_{\Xtrunc}(\Jac) -2d  - 2g \\
&=& -1,
\end{eqnarray*}
which is generically empty. Thus, it follows that $[A_N,\Phi_N]$ must
be reducible. Moreover, according to
Proposition~\ref{cor:DimNSLocus}, 
\begin{eqnarray*}
\edim
\Restrict(\Theta) \cap
\Restrict(\ModSW_{\Xtrunc}(\Jac)\cap\Divisors(a))&=& 
2\mult + \edim \ModSW_{\Xtrunc}(\Jac) -2d  - 2g \\
&=& -2,
\end{eqnarray*}
which is also generically empty. Hence, $[A_N,\Phi_N]$ and
$[A_{\Xtrunc},\Phi_{\Xtrunc}]$ 
satisfy Hypotheses~(H-\ref{item:Red}).
\item
\label{item:IrrIrr}
If $\Restrict[A_N,\Phi_N]$ and 
$\Restrict[A_{\Xtrunc},\Phi_{\Xtrunc}]$ both lie in
$\CritMan$, then the Hypotheses~(H-\ref{item:Irr}) are satisfied.
\end{list}

The assertion at the beginning of the proposition follows easily.
\end{proof}

The above proposition says that we can partition the points in the
cut-down moduli space (which is an oriented, zero-dimensional
manifold) for sufficiently large $T$ into two disjoint sets, the
subsets of configurations which satisfy  (H-\ref{item:Red}) and
(H-\ref{item:Irr}) respectively. Thus, if we let $\SWred_{\spinc}(a)$ and
$\SWirr_{\spinc}(a)$ be the signed number of points satisfying (H-1)
and (H-2) respectively, then
\begin{equation}
\label{eq:RedIrr}
\SW_{X,\spinc}(a)=\SWred_{\spinc}(a)+\SWirr_{\spinc}(a).
\end{equation} As we
shall see, gluing theory allows us to compute both of these quantities
in terms of cylindrical-end moduli spaces. So, in the next step, we
study these cylindrical-end moduli spaces.

\begin{lemma}
For all $\SpinC$ structures $\spinc$ on $X$
the corresponding moduli spaces $\ModSW_N(\CritMan)$, $\ModSW_{\Xtrunc}(\Jac)$, and
$\ModSW_{\Xtrunc}(\CritMan)\cap \Divisors(a)$ are all compact manifolds.
\end{lemma}

\begin{proof}
The compactness of $\ModSW_{\Xtrunc}(\Jac)$ and $\ModSW_N(\CritMan)$
follows from the usual compactness arguments, together with the facts 
that the Chern-Simons-Dirac functional is real-valued,
$\CSD(\Jac)>\CSD(\CritMan)$, and there are no other critical
manifolds. Compactness of $\ModSW_{\Xtrunc}(\CritMan)\cap
\Divisors(a)$ follows from this, together with a straightforward dimension count
(see the discussion above in the proof of Lemma~\ref{lemma:TwoCases},
part (P-\ref{item:RedIrr})).
\end{proof}

Compactness of $\ModSW_{\Xtrunc}(\Jac)$ allows us to define a relative
invariant with reducible boundary values. We pause to discuss some
relevant properties of this invariant.

\begin{defn}
Let $\spinc_0$ be a $\SpinC$ structure on $\Xtrunc$ which extends over
$X$. 
Since the moduli space $\ModSW_{\Xtrunc,\spinc_0}(\Jac)$ is compact,
there is a {\em relative Seiberg-Witten invariant}
$$\SW_{(\Xtrunc,\spinc_0,\Jac)}\colon \Alg(\Xtrunc) \longrightarrow \Z,$$
defined by the pairing
$\SW_{(\Xtrunc,\spinc_0,\Jac)}(a)=\langle
[\ModSW_{\Xtrunc,\spinc_0}(\Jac)],\mu(a)\rangle.$
\end{defn}

This relative invariant is related to an absolute invariant, according
to the following.

\begin{prop}
\label{prop:SmallSpinc}
If $\spinc$ satisfies 
$-n< \langle c_1(\spinc),[\Sigma]\rangle \leq n$, then for all $a\in\Alg(X)$,
$$\SW_{X,\spinc}(a)=\SW_{(\Xtrunc,\spinc_0,\Jac)}(a),$$
where $\spinc_0=\spinc|_{\Xtrunc}$.
\end{prop}

\begin{proof}
Recall that $\ModSW(\spinc|_{N})$ consists entirely of reducibles all of
which are smooth, according to Proposition~\ref{prop:UnobstructedOnN};
thus, gluing theory identifies the moduli spaces
$\ModSW_{X(T)}(\spinc)$ for large $T$ with
$\ModSW_{\Xtrunc,\spinc_0}(\Jac)$. (See
also~\cite{SympThom}, where this result appears as Proposition~2.7.)
\end{proof}

We now return to the discussion of $\SWred$ and $\SWirr$.  Although
the definitions of both terms implicitly use $T$, we show now that if
$T$ is sufficiently large, then the terms can be computed from
absolute invariants (and hence are independent of the parameter).

\begin{prop}
\label{prop:RedBoundaryTwo}
Suppose that $\spinc$ satisfies $$\BigSpinc,$$
where $\Sigma$ has self-intersection number $-n$,
and let
$$e=g-1+\frac{\langle c_1(\spinc),[\Sigma]\rangle + n}{2}.$$
Then, for all sufficiently large $T$, 
$$\SWred_{\spinc}(a) =
SW_{X,\spinc-\PD[\Sigma]}(a \cm \xi_{g-1-e}(\Sigma)),$$
where $\xi_{g-1-e}(\Sigma)\in\Alg(\Sigma)$ is the element defined in
Definition~\ref{def:xis}. 
\end{prop}

\begin{proof}
The moduli space $\ModRed_N(\Jac)-\Theta$ comes equipped with an
obstruction bundle $\OBundle \longrightarrow \ModRed_N(\Jac)-\Theta$, defined
by $\OBundle_{[(A,0)]}=\CoKer \Dirac_A$.  (whose $K$-theory class canonically
extends over all of $\ModRed_N(\Jac)$).  The dimension count in
Lemma~\ref{lemma:TwoCases} guarantees that each solution in
$\ModSW_{\Xtrunc}(\Jac)\cap \Divisors(a)$ extends uniquely to a smooth
reducible over $N$. Thus, gluing theory gives that
\begin{eqnarray*}
\SWred_{\spinc}(a)
&=&
\langle [\ModSW_{\Xtrunc}(\Jac)\cap \Divisors(a)], \Euler({\mathcal L}\otimes \Restrict^*(\OBundle)) \rangle \\
&=& \langle [\ModSW_{\Xtrunc}(\Jac)], \mu(a)\cup \Euler({\mathcal L}\otimes
\Restrict^*(\OBundle))\rangle,
\end{eqnarray*}
where $\Euler$ denotes the Euler class of a bundle (or $K$-theory
element).The Riemann-Roch formula says
$\dim(\OBundle)=2g-2-2e$.
Using the index theorem for families, together with the
holomorphic interpretation of the obstruction bundle $\OBundle$ given in
Proposition~\ref{prop:KerDirac}, it is a straightforward computation that the
total Chern class of $\OBundle$ is
$$\prod_{i=1}^g(1+\mu(A_i)\mu(B_i))$$
(see also~\cite{SympThom} Proposition~2.6);
thus, 
$$\Euler({\mathcal L}\otimes \Restrict^*(\OBundle))=
c_{g-1-e}({\mathcal L}\otimes \Restrict^*(\OBundle))=
\xi_{g-1-e}([\Sigma]).$$
Putting all this together, we have that
\begin{equation}
\label{eq:ObstructedVersion}
\SWred_{\spinc}(a)=\SW_{(\Xtrunc,\spinc_0,\Jac)}(a \cm \xi_{g-1-e}(\Sigma)),
\end{equation}
where $\spinc_0=\spinc|_{\Xtrunc}$. 
Since $n-2g+2\leq \langle c_1(\spinc-\PD[\Sigma]), [\Sigma]\rangle
\leq n$ and 
$\spinc-\PD[\Sigma]|_{\Xtrunc}=\spinc_0$,
the proposition then follows from Proposition~\ref{prop:SmallSpinc}.
\end{proof}

\begin{prop}
\label{prop:SWirrIsRelativeInvariant}
For sufficiently large $T$,
$$
\SWirr_{\spinc}(a)
=
\# \ModSW_{\Xtrunc}(\CritMan)\cap \Divisors(a).
$$
\end{prop}

\begin{proof}
Gluing shows that
$$
\SWirr_{\spinc}(a)
=
\Big(\# \ModSW_{\Xtrunc}(\CritMan)\cap \Divisors(a)\Big) \Big(\deg(\Restrict\colon
\ModSW_{\NSig}(\CritMan)\rightarrow \CritMan)\Big).
$$
According to Propositions~\ref{prop:NSigDefTheory} and
\ref{prop:DimensionCountingPert}, 
$\Restrict\colon
\ModSW_{\NSig}(\CritMan)\rightarrow \CritMan$
either has degree $+1$, or $\ModSW_{\NSig}(\CritMan)$ is empty. The
latter case would force $\SWirr_{\spinc}(a)\equiv 0$ (for the given
genus and self-intersection number).

To rule out this latter case, we need only look at an example where
the irreducible term is non-zero.  Let $X$ be a ruled surface $X$ over
$\Sigma$ associated to the line bundle with Euler number $-n$. Let
$\Sigma\subset X$ denote the section with self-intersection number
$-n$, and fix any $0\leq e \leq g-1$.  Let $\spinc$ denote the
$\SpinC$ structure over $X$ given by $\spinc=\spinc_0 + e\PD[F]$,
where $\spinc_0$ is the canonical $\SpinC$ structure on $X$ associated
to the K\"ahler structure, and $F$ denotes a fiber in the ruling. It
is easy to see that $\SW_{X,\spinc-\PD[\Sigma]}\equiv 0 $, as the
corresponding space of divisors is empty (see
Proposition~\ref{prop:RSurfCohomology}). Moreover, we know that
$\SW_{X,\spinc}\not\equiv 0$ (compare
Example~\ref{subsec:RSurfEx}). Thus, in light of
Equation~\eqref{eq:RedIrr} and Proposition~\ref{prop:RedBoundaryTwo},
we have examples where $\SWirr\not\equiv 0$, forcing the degree to be
non-zero.
\end{proof}

We will give the seemingly {\em ad hoc} quantity $\#
\ModSW_{\Xtrunc}(\CritMan)\cap \Divisors(a) $ a more intrinsic formulation in
Section~\ref{sec:SWirr}. With the help of this formulation, we can
then prove a vanishing result for this term under suitable
algebro-topological hypotheses on the embedding of $\Sigma\subset X$
(Proposition~\ref{prop:Vanishing}). 

\section{Relative Invariants}
\label{sec:SWirr}

Let $\Xtrunc$ be a smooth, oriented manifold-with-boundary with
$b_2^+(\Xtrunc)>0$,
whose boundary is identified with $\partial \Xtrunc
= -Y$, a circle bundle over a Riemann surface $\Sigma$ of genus $g>0$
with Euler number $-n$, where
$n> 2g-2$.

In Section~\ref{sec:ProductFormula}, we studied the moduli space
$\ModSW_{\Xtrunc}(\CritMan)$, and used it to define a relative
invariant $$\SWirr_{\spinc}\colon\Alg(\Xtrunc)\longrightarrow \Z,$$ by
cutting down the moduli space $\ModSW_{\Xtrunc}(\CritMan)$ by
submanifolds representing $\mu(a)$ which are induced from compactly
supported representatives for homology in $\Xtrunc$ (see
Proposition~\ref{prop:SWirrIsRelativeInvariant}). When $a\in\Alg(Y)$,
there are alternate representatives which are supported ``at
infinity.'' The advantage of these representatives is that the
corresponding relative invariant inherits relations arising from the
cohomology ring of $\CritMan$.  In view of the non-compactness of
$\ModSW_{\Xtrunc}(\CritMan)$, the two types of representatives do not
necessarily give rise to the same invariant. However, the difference
can be explicitly computed in terms of other Seiberg-Witten
invariants. In this section, we recast this discussion in a more
algebraic setting, defining an invariant
$$\SW_{(\Xtrunc,\CritMan)}\colon
\Alg(\Xtrunc)\otimes H^*(\CritMan) \longrightarrow
\Z.$$
which simultaneously captures both types of representatives; in
particular, $$\SWirr_{\spinc}(a)=\SW_{(\Xtrunc,\CritMan)}(a\otimes 1).$$
Proposition~\ref{prop:Vanishing} then follows from properties of this
invariant.

A subtlety arises in the definition of $\SW_{(\Xtrunc,\CritMan)}$, since
the moduli space $\ModSW_{\Xtrunc}(\CritMan)$ is not compact. However, we have
the following weak compactness theorem.

\begin{defn}
A sequence of configurations $\{[A_i,\Phi_i]\}_{i=1}^{\infty}$ is said
to {\em converge weakly} to a configuration
$$[B,\Psi]\times[A,\Phi]\in 
\UnparModFlow(\CritMan,\Jac)\times_{\Jac}\ModSW_{\Xtrunc}(\Jac)$$
if $[A_i,\Phi_i]$ converges to $[A,\Phi]$ in $\CinftyLoc$, and there
is an increasing, unbounded sequence of real numbers
$\{T_i\}_{i=1}^{\infty}$ with $T_i>i$, so that the translates of
$\{[A_i,\Phi_i]|_{[0, 2T_i]\times Y}\}_{i=1}^{\infty}$, viewed as a
sequence of configurations on $[-T_i,T_i]\times Y$, converge in
$\CinftyLoc$ to a configuration which is equivalent (under
translations) to $[B,\Psi]$.
\end{defn}

\begin{prop}
\label{prop:WeakCompactness}
Weak convergence gives the space
$$\CompactModSW_{\Xtrunc}(\CritMan)=\ModSW_{\Xtrunc}(\CritMan)\coprod
\UnparModFlow(\CritMan,\Jac)\times_{\Jac}\ModSW_{\Xtrunc}(\Jac)$$
the structure of a compact Hausdorff space.
\end{prop}

\begin{proof}
This a standard argument from Morse-Floer theory. A general discussion
of compactness results for the anti-self-duality equation can be found
in~\cite{MMR} (see especially Theorem~6.3.3 of~\cite{MMR}); so we
sketch the argument here only briefly.

A sequence $[A_i,\Phi_i]\in\ModSW_{\Xtrunc}(\CritMan)$ converges
in $\CinftyLoc$, after passing to a subsequence, to some solution
$[A,\Phi]$ to the Seiberg-Witten equations on $\Xtrunc$. Since each of
the $[A_i,\Phi_i]$ have finite energy, so does $[A,\Phi]$; thus, it
has a boundary value. If $\Restrict[A,\Phi]\in\CritMan$, then the
length-energy estimates of L. Simon (~\cite{Simon})
can be used to show that the convergeance is $\Cinfty$ as in \cite{MMR}.

If, on the other hand, $\Restrict[A,\Phi]\not\in\CritMan$, it must be
the case that $\Restrict[A,\Phi]\in\Jac$. Now, let $T_i\in\R$ be the
number so that $$\CSD[A_i,\Phi_i]_{\{T_i\}\times
Y}=\frac{\CSD(\Jac)+\CSD(\CritMan)}{2}.$$ Clearly, $T_i\goesto
\infty$.  After passing to a subsequence, we can find a configuration
$[B,\Psi]$ so that the sequence $[A_i,\Phi_i]|_{[0,2T_i]\times Y}$,
viewed as a sequence of configurations over $[-T_i,T_i]$, converges in
$\CinftyLoc$ to $[B,\Psi]$. In fact, $[B,\Psi]$ must solve the
Seiberg-Witten equations and it must have finite energy, so
$[B,\Psi]\in\ModFlow(\CritMan,\Jac)$.  The usual length-energy
estimates then guarantee that the boundary values match up.
\end{proof}

The topological space $\CompactModSW_{\Xtrunc}(\CritMan)$ defined in
Proposition~\ref{prop:WeakCompactness} is called the compactified
moduli space. The following result follows immediately from its
definition. 

\begin{prop}
The inclusion maps 
$$\Include\colon\ModSW_{\Xtrunc}(\CritMan)\longrightarrow\SWConfigIrr(\Xtrunc-(0,\infty)\times
Y)$$
and 
$$\Include\circ\Proj_2 \colon
\UnparModFlow(\CritMan,\Jac)\times_{\Jac}\ModSW_{\Xtrunc}(\Jac) \longrightarrow
\SWConfigIrr(\Xtrunc-(0,\infty)\times Y)$$
fit together to give a continuous map $$\CompactInclude\colon
\CompactModSW_{\Xtrunc}(\CritMan)\longrightarrow \SWConfigIrr(\Xtrunc-(0,\infty)\times
Y),$$
where $\SWConfigIrr$ denotes the irreducible configurations.
\end{prop}

Similarly, we can extend the restriction map over the compactified
moduli space, as follows.

\begin{prop}
The restriction maps 
$$\Restrict_{\CritMan}\colon\ModSW_{\Xtrunc}(\CritMan)\longrightarrow
\CritMan$$
and 
$$\Restrict_{\CritMan}\circ\Proj_1 \colon
\UnparModFlow(\CritMan,\Jac)\times_{\Jac}\ModSW_{\Xtrunc}(\Jac) \longrightarrow
\CritMan$$
fit together to give a continuous map $$\CompactRestrict_{\CritMan}\colon
\CompactModSW_{\Xtrunc}(\CritMan)\longrightarrow \CritMan.$$
\end{prop}

\begin{proof}
If a sequence
$[A_n,\Phi_n]\in\ModSW_{\Xtrunc}(\CritMan)$ converges to an ideal
point
$$[B,\Psi]\times[A,\Phi]\in\UnparModFlow(\CritMan,\Jac)\times_{\Jac}\ModSW_{\Xtrunc}(\Jac),$$
then there is a divergent sequence $\{T_n\}_{n=1}^{\infty}$ of real
numbers so that
$$\lim_{n\goesto\infty}\tau_n^*[A_n,\Phi_n]|_{[T_n,\infty)\times Y} =
[B,\Psi]|_{[0,\infty)\times Y},$$ where $$\tau_n\colon
[0,\infty)\times Y \longrightarrow [T_n,\infty)\times Y$$ is the map
induced by translation by $T_n$ on the first coordinate.  Since each
path has finite energy, continuity of the restriction maps
(see~\cite{MMR}) guarantees
that $$\lim_{n\goesto\infty}\Restrict[A_n,\Phi_n]|_{\{t\}\times Y} =
\lim_{n\goesto\infty}\Restrict\tau_n^*[A_n,\Phi_n]=
\Restrict[B,\Psi].$$
\end{proof}

Gluing gives this space more structure. 

\begin{prop}
\label{prop:MovePointOne}
Gluing endows $\CompactModSW_{\Xtrunc}(\CritMan)$ with the structure
of a manifold. The space of ideal solutions
$$\UnparModFlow(\CritMan,\Jac)\times_{\Jac}\ModSW_{\Xtrunc}(\Jac)$$
has the structure of a smooth submanifold of codimension two.  In
particular, a fundamental class for $\ModSW_{\Xtrunc}(\CritMan)$ gives
rise to a unique fundamental class for
$\CompactModSW_{\Xtrunc}(\CritMan)$.
\end{prop}

\begin{proof}
Gluing describes the end of $\ModSW(X,\CritMan)$ as a fibered product
$$\Big(\UnparModFlowBased(\CritMan,\Jac)\times_{\Jac}
\ModSWBased_{\Xtrunc}(\Jac)\times (0,\infty)\Big)/S^1,$$
where the superscript denotes based versions of the moduli
spaces. This gives the space of ideal solutions a disk-bundle
neighborhood in $\CompactModSW_{\Xtrunc}(\CritMan)$.
\end{proof}

In light of the above result, we can define the relative Seiberg-Witten
invariant $\SW_{(\Xtrunc,\CritMan)}$, as follows.

\begin{defn}
The {\em relative Seiberg-Witten invariant}
$$\SW_{(\Xtrunc,\CritMan)}\colon \Alg(\Xtrunc)\otimes H^*(\CritMan) \longrightarrow
\Z$$
is defined by
$$\SW_{(\Xtrunc,\CritMan)}(a\otimes\omega)=\langle
[\CompactModSW_{\Xtrunc}(\CritMan)],
\CompactInclude^*(\mu(a))\cup \CompactRestrict_{\CritMan}^*(\omega)\rangle.$$
\end{defn}

We now spell out the strategy for proving
Proposition~\ref{prop:Vanishing}.  First, it is shown that for
$b\in\Alg(Y)$, $\SW_{(\Xtrunc,\CritMan)}(a\cm b\otimes
\omega)$ can be expressed in terms of $\SW_{(\Xtrunc,\CritMan)}(a\otimes
b\cm\omega)$ and $\SW_{(\Xtrunc,\Jac)}$
(Lemma~\ref{lemma:Commutator} and
Proposition~\ref{prop:MovePoint}). Here, $b\cm\omega$ denotes the
action of $\Alg(Y)$ on $H^*(\CritMan)$ induced from the inclusion of
$\CritMan$ in $\SWConfigIrr(Y)$.  (Note the cohomology classes over
$\CompactModSW_{\Xtrunc}(\CritMan)$ induced from $\Alg(Y)$ through the
action on $H^*(\CritMan)$, and pulled back via $\CompactRestrict$,
correspond to divisor representatives over $\Xtrunc$ which are
supported ``at infinity.'')  Then, it is shown that
$\SW_{(\Xtrunc,\CritMan)}(a \otimes b\cm\omega)$ vanishes, when $b$ has
sufficiently high degree. This follows from algebraic considerations,
according to which $b\cm\omega=b'\cm\omega$, where $b'\in\Alg(Y)$ lies
in the ideal generated by the cycles in $Y$ which bound in $\Xtrunc$
(Proposition~\ref{prop:Algebra}). It is then easy to see that
$\SW(a \otimes b'\cm\omega)$ vanishes (Corollary~\ref{cor:VanishingOfRelInv}).

Now, we express the ``commutator'' $\SW(a \otimes b\cm\omega)-\SW(a\cm
b \otimes \omega)$. First note that if $b$ is induced from $H_1(Y)$,
the commutator vanishes, as follows.

\begin{lemma}
\label{lemma:Commutator}
Let $[\gamma]\in H_1(Y)$, then for all $a\in\Alg(X)$ and $\omega\in
H^*(\CritMan)$, 
$$\SW_{(\Xtrunc,\CritMan)}(a\otimes \mu[\gamma]\cm\omega)=
  \SW_{(\Xtrunc,\CritMan)}(a \cm\mu[\gamma] \otimes\omega).$$
\end{lemma}

\begin{proof} We must show that $\CompactRestrict_{\CritMan}^*(\mu[\gamma])$ is
homologous to $\CompactInclude^*(\mu[\gamma])$. It suffices to
verify this over the subset $\ModSW(\Xtrunc,\CritMan)\subset
\CompactModSW_{\Xtrunc}(\CritMan)$, since the complement has
codimension two, and the classes in question are one-dimensional. 
Over the subset, now, the claim is easy to verify. On
$\ModSW_{\Xtrunc}(\CritMan)$, $\Restrict^*(\mu[\gamma])$ is
represented by $(\Hol_\gamma\circ
\Restrict_{\CritMan})^*(d\theta)$, the holonomy around a representative of
$\gamma$ ``at infinity'' (see Proposition~\ref{prop:MuCircle});
while $\Include^*(\mu[\gamma])$ is represented by
$\Hol_{\gamma_0}^*(d\theta)$,
where $\gamma_0=0\times\gamma\subset 0\times Y\subset [0,\infty)\times
Y\subset X$. Now, the cylinder $[0,\infty)\times \gamma$ provides a
homotopy between $\Hol_\gamma\circ \Restrict_{\CritMan}$ and $\Hol_{\gamma_0}$.
\end{proof}

It remains to see how the point class commutes. For this class, we can
express the commutator in terms of $\SW_{(X,\Jac)}$ and
another Seiberg-Witten invariant, 
defined below.

\begin{defn}
There is a {\em Seiberg-Witten invariant of the tube}
$$\SWTube\colon H^*(\CritMan)\longrightarrow H^*(\Jac)\subset
\Alg(\Xtrunc)$$
which raises degree by dimension $2g-\dim\UnparModFlow(\CritMan,\Jac)=
2g-\dim\CritMan$, defined by
$$\SWTube(\omega)=(P_2)_*\Big((\Restrict_{\Jac}\times\Id)^*\PD[\Delta]\cup
(\Restrict_{\CritMan}\circ P_1)^*\omega\Big),$$
where $P_1$ and $P_2$ are the projection maps
\begin{eqnarray*}
P_1\colon \UnparModFlow(\CritMan,\Jac)\times\Jac
\longrightarrow \UnparModFlow(\CritMan,\Jac) &{\text{and}}&
P_2\colon \UnparModFlow(\CritMan,\Jac)\times\Jac
\longrightarrow \Jac,
\end{eqnarray*}
and $\PD[\Delta]$ denotes the Poincar\'e dual of the diagonal
$\Delta\subset \Jac\times\Jac$.
Thus, $\SWTube$ satisfies:
$$\langle\UnparModFlow(\CritMan,\Jac)\times_{\Jac}\ModSW_{\Xtrunc}(\Jac),
\CompactInclude^*(\mu(a))\cup\CompactRestrict_{\CritMan}^*(\omega)\rangle =
\SW_{(\Xtrunc,\Jac)}(a\cm\SWTube(\omega)).$$
\end{defn}

We can now calculate the commutator, which involves comparing the
cohomology classes $\CompactRestrict_{\CritMan}^*\mu(y)$ and
$\CompactInclude^*\mu(x)$ over $\CompactModSW_{\Xtrunc}(\CritMan)$,
where $y$ is a point in $Y$ and $x$ is a point in $\Xtrunc$.

\begin{prop}
\label{prop:MovePoint}
Choose points $x\in \Xtrunc$ and $y\in Y$.
In $\CompactModSW_{\Xtrunc}(\CritMan)$, we have
$$\CompactRestrict_{\CritMan}^*(\mu(y))-\CompactInclude^*(\mu(x)) =
\PD[\UnparModFlow(\CritMan,\Jac)\times_{\Jac}\ModSW_{\Xtrunc}(\Jac)].$$
Consequently, there is a relation between Seiberg-Witten invariants:
$$\SW_{(\Xtrunc,\CritMan)}(a\otimes\mu(y)\cm\omega)-\SW_{(\Xtrunc,\CritMan)}(a\cm\mu(x)
\otimes \omega) = \SW_{(\Xtrunc,\Jac)}(a\cm \SWTube(\omega)).$$
\end{prop}

\begin{proof}
Clearly, the difference
$\CompactRestrict_{\CritMan}^*(\mu(y))-\CompactInclude^*(\mu(x))$ is
the first Chern class of the circle bundle
$\Hom_{S^1}({\BaseFib_x},{\BaseFib_y})$.  Here, $\BaseFib_z$ denotes
the moduli space based at $z$; see Section~\ref{sec:Cohomology}. To
prove the proposition, we must verify that this bundle admits a
section $\sigma$ in the complement of
$$\UnparModFlow(\CritMan,\Jac)\times_{\Jac}\ModSW_{\Xtrunc}(\Jac)\subset
\CompactModSW_{\Xtrunc}(\CritMan)$$
(i.e. over $\ModSW_{\Xtrunc}(\CritMan)\subset
\CompactModSW_{\Xtrunc}(\CritMan)$) and that, with respect to a
trivialization of the circle bundle over a disk transverse to the
submanifold, the restriction of the section to the boundary induces a
map from the circle to the circle which has degree one.

The section $\sigma$ is induced by parallel transport, as
follows. Let $\gamma$ be a half-infinite arc
formed by joining $[0,\infty)\times y$ to
any arc which connects 
$x$ to $0\times y$. Over the point $[A,\Phi]\in\ModSW_{\Xtrunc}(\CritMan)$,
parallel transport via $A$ along $\gamma$ induces a homomorphism 
in $\Hom_{S^1}({\BaseFib_x},{\BaseFib_y})$. 

We now verify that the trivialization induces a degree one map around
circles transverse to the submanifold. For any point in the submanifold
$$[A_1,\Phi_1]\times[A_2,\Phi_2]\in
\UnparModFlow(\CritMan,\Jac)\times_{\Jac}\ModSW_{\Xtrunc}(\Jac),$$
fix fibers
\begin{eqnarray*}
[A_1,\Phi_1,\lambda_1]\in \BaseFib_x|_{[A_1,\Phi_1]}
&{\text{and}}&
[A_2,\Phi_2,\lambda_2]\in \BaseFib_y|_{[A_2,\Phi_2]}.
\end{eqnarray*}
These choices induce a trivialization of
$\Hom_{S^1}({\BaseFib_x},{\BaseFib_y})$ over a disk in
$\CompactModSW_{\Xtrunc}(\CritMan)$ transverse to
$[A_1,\Phi_1]\times[A_2,\Phi_2]$ (obtained by varying the gluing and
translation parameters). Calculating the desired degree amounts to
seeing how the holonomy along $\gamma$ varies as the gluing parameter
is rotated. But holonomy along any path which crosses the gluing
region once varies as a degree one function of the gluing parameter.
\end{proof}

We can understand the action of $\Alg(Y)$ on $H^*(\CritMan)$
explicitly, under the identification $\CritMan\cong \Sym^k(\Sigma)$.

Before describing this, we begin with a few preliminaries about the homology
of symmetric products of $\Sigma$ (for an extensive discussion of this
topic, see~\cite{MacDonald}). Recall that $\Sym^k(\Sigma)$ can be
viewed as the quotient of the $k$-fold Cartesian product
$\Sigma^{\times k}$ by the action of
the symmetric group on $k$ letters. We denote the quotient map by
$$\QuotMap\colon \Sigma^{\times k}
\longrightarrow \Sym^k(\Sigma).$$ According to elementary properties
of the transfer homomorphism, $$\QuotMap_*\colon H_*(\Sigma^{\times k})\longrightarrow
H_*(\Sym^k(\Sigma))$$ is surjective.  Dually, we have a 
map $$\QuotMap^*\colon H^*(\Sym^k(\Sigma))\longrightarrow
H^*(\Sigma^{\times k})$$ which identifies $H^*(\Sym^k(\Sigma))$ with
the elements of $H^*(\Sigma^{\times k})\cong H^*(\Sigma)^{\otimes k}$
which are invariant under the symmetric group action. In particular,
by summing over the action, we obtain a map $$\Sym^k\colon
H^*(\Sigma)\longrightarrow H^*(\Sym^k(\Sigma)).$$ 
Thus, if we fix any  collection of points $\{p_2,...,p_k\}\subset
\Sigma$, given
$\omega\in H^*(\Sigma)$, $\Sym^k(\Sigma)$ is the class characterized
by the property that 
$$
\langle \Sym^k(\omega), \QuotMap_*(Z\times p_2 \times
... \times p_k)\rangle  = \langle \omega,Z\rangle,$$ for any cycle $Z\subset
\Sigma$.  Equivalently, given a cycle $Z\in H_*(\Sigma)$, $\Sym^k(\PD[Z])$ is
Poincar\'e dual to the cycle $\QuotMap(Z\times \Sigma\times ...\times
\Sigma)$.
The above discussion works over rational coefficients (which suffices
for our purposes), but in fact it works over $\Z$ as well, since
$H_*(\Sym^k(\Sigma))$ has no torsion (see~\cite{MacDonald}).

\begin{prop}
\label{prop:MuClassesOnSym}
Under the identification $\CritMan\cong \Sym^k(\Sigma)$, the canonical
map $$\Alg(Y)\longrightarrow H^*(\Sym^k(\Sigma)),$$ induced from the
inclusion of $\Sym^k(\Sigma)=\CritMan\longrightarrow \SWConfigIrr(Y)$,
takes $\mu(y)$ for $y\in H_*(Y)$ to the cohomology class
$\Sym^k(\PD[\pi_*(y)])\in H^{2-*}(\Sym^k(\Sigma))$, where $\pi\colon
Y\rightarrow \Sigma$ is the projection map.
\end{prop}

\begin{proof}
We can reduce to a corresponding statement for configurations over
$\Sigma$, as follows. Let $E$ be a line bundle over $\Sigma$, so that 
$W\cong \pi^*(E\otimes (\C\oplus\Canon{\Sigma}))$. Then, pull-back induces a
map
$$\pi^*\colon \SWConfig(\Sigma,E)=\Conns(E)\times
\Sections(E)/\Map(\Sigma,S^1)\longrightarrow \SWConfig(Y,W),$$
to the configurations where the fiber-wise holonomy of the connection
is constant, and the section is covariantly constant around each
fiber. The identification between the critical manifolds and the
symmetric powers $\CritMan\cong \Sym^k(\Sigma)$ described
in~\cite{MOY} is obtained by proving that $\CritMan$ lies in the image
of this pull-back map, and indeed that it lies in the pull-back of the
vortex moduli space, which, according to~\cite{JaffeTaubes} (see
also~\cite{Bradlow}), is in turn identified with the space of
divisors, by looking at the zero-set of the section.  The key points
we need presently are that $\CritMan$ lies in
$\pi^*(\SWConfig(\Sigma,E))$, and that configurations are the
pull-backs of configurations $[A,\Phi]\in\SWConfig(\Sigma,E)$, where $\Phi$
is $\DBar_A$-holomorphic section.

Over $\SWConfig(\Sigma,E)$, there is a universal line bundle ${\mathcal
L}(\Sigma)$, defined in the usual manner. Note that 
$${\mathcal L}(Y)|_{\pi^*(\SWConfig(\Sigma))\times Y}
\cong \pi^*({\mathcal L}(\Sigma)),$$
so $\mu(y)|_{\pi^*(\SWConfig(\Sigma))}$ for
$y\in H_*(Y)$ agrees with $\mu([\pi_*(y)])$, where the former
$\mu$-map is induced from ${\mathcal L}(Y)$, and the latter from
${\mathcal L}(\Sigma)$.  We have thus reduced the proof of the
proposition to a statement purely over $\Sigma$; so for the duration
of the proof, ${\mathcal L}$ will refer to ${\mathcal L}(\Sigma)$,
$\SWConfig$ will refer to $\SWConfig(\Sigma,E)$, and all $\mu$-maps
will be calculated over $\Sigma$.

To facilitate the proof over $\Sigma$, we pause for a discussion about
the canonical section $\sigma$ of the universal line bundle ${\mathcal
L}$, which takes the configuration
$[A,\Phi]\times\{x\}\in\SWConfig\times \Sigma$ to the based
configuration $[A,\Phi,\Phi(x)]$. This section has the property that,
under the canonical identification of ${\mathcal L}|_{[A,\Phi]\times
\Sigma}\cong E$ (where $E$ is the bundle over $\Sigma$ with Chern
number $k$), the restriction $\sigma|_{[A,\Phi]\times\Sigma}$ is
identified with the section $\Phi$ of $E$. In particular, if
$\Phi\not\equiv 0$ is a holomorphic section, then
$\sigma|_{[A,\Phi]\times \Sigma}$ has at most $k$ zeros; moreover, if
it has $k$ zeros, then each is transverse.

Now, if $y\in \Sigma$ is a point (i.e. a generator of $H_0(\Sigma;\Z)$),
then by definition, $\mu(y)$ is the element of $H^2(\CritMan)$ whose
pairing against any homology class $[S]\in H_2(\CritMan)$ is given by
$$\langle \mu(y), [S]\rangle =
\langle c_1({\mathcal L}_y),
[S]\rangle.$$ where, as usual, ${\mathcal L}_y$ denotes the
restriction of
${\mathcal L}$ to ${\SWConfig(\Sigma,E)}\times \{y\}$. Choose points
$\{p_2,...,p_k\}\subset \Sigma$ which are distinct from $y$.
Recall that $H_2(\Sym^k(\Sigma))$ is generated by the
surface $\QuotMap(\Sigma\times p_2\times ...\times p_k)$ (where
$p_2,...,p_k$ are points on $\Sigma$), and the tori of the form
$\QuotMap(C_1\times C_2
\times p_3\times ...\times p_k)$, where $C_1,C_2\subset \Sigma$ are closed curves
in $\Sigma$, which we can choose to miss $y$. 
The canonical section $\sigma$ restricted to  a torus of the form 
$\QuotMap(C_1\times C_2
\times p_3\times ...\times p_k)\times\{y\}$ clearly vanishes nowhere (as all $k$ of the zeros
have been constrained to lie in the set $C_1\cup C_2
\cup\{p_3,...,p_k\}$ which does not include the point $y$); thus,
$$\langle c_1({\mathcal L}_y), [\QuotMap(C_1\times C_2 \times...\times
p_k)]\rangle = 0.$$

Over $\QuotMap(\Sigma\times p_2,...\times p_k)$ the canonical section
vanishes at the single point $\QuotMap(y\times p_2\times...\times
p_k)$. We verify transversality of this zero, as follows.  View
$\sigma$ as a section over $\Sigma\times\Sigma=q(\Sigma\times
p_2\times ... \times p_k)\times \Sigma$; we know that
$\sigma(y,y)\equiv 0$, and that $\Deriv\sigma_{(y,y)}$ induces an
isomorphism from $0\oplus T_y\Sigma$ to $E_y$ (i.e. that the zero of
$\sigma|_{\{y\}\times\Sigma}$ at $y$ is transverse). Differentiating
the equation that $\sigma(y,y)\equiv 0$, we see that
$$\Deriv\sigma_{(y,y)}(0,v)=-\Deriv\sigma_{(y,y)}(v,0).$$ Thus, the
section $\sigma|_{\Sigma \times \{y\}}$ of ${\mathcal
L}_y|_{[\Sigma]}$ has a single, transverse zero, which shows that
$$\langle c_1({\mathcal L}), [\QuotMap(\Sigma\times p_2\times...\times
p_k)]\rangle = \pm 1. $$ Moreover, the sign is positive since the
section is holomorphic.

Hence, we have proved the result when $y\in H_0(\Sigma)$. Proving the
result for classes coming from $H_1(\Sigma)$ amounts to proving that,
if $C_1$ and $C_2$ closed curves in $\Sigma$ which meet transversally,
then $$\langle c_1({\mathcal L}),\QuotMap(C_1\times p_2 \times ... \times
p_n)\times C_2\rangle = -\# C_1\cap C_2
= \# C_2 \cap C_1.$$

Note first that the zeros of the canonical section $\sigma$,
restricted to $\QuotMap(C_1\times p_2 \times ... \times p_k)\times C_2$
are the points $C_1\cap C_2$ (a zero of $\sigma$ corresponds to a
point where the section $\Phi$ vanishes at some point of $C_2$, but
the zeros of $\Phi$ lie in $C_1\cup \{p_2,...,p_k\}$, and
$\{p_2,...,p_k\}\cap C_2$ is empty). We must now consider the local
contribution of each zero (and check transversality).

Consider the map 
$C_1\times C_2 \colon S^1\times S^1\longrightarrow
\Sym^k(\Sigma)\times \Sigma$
defined by $C_1\times C_2(s,t)=\QuotMap(C_1(s)\times p_2\times...\times
p_k)\times C_2(t)$. Suppose for notational simplicity that
$C_1(0)=C_2(0)=y$.  We can view $\sigma$ as a section of ${\mathcal L}$
pulled back to this torus. Now, evaluated on a typical tangent vector
to the torus $a\DDs+b\DDt$, the derivative of $\sigma$ at the
intersection point is given by
\begin{eqnarray}
\Deriv_{(y,y)}\sigma\circ (C_1\times C_2) (a\DDs + b\DDt)
&=&
a \Deriv_{(y,y)} \sigma(\dd{C_1}{s}(0),0) + 
b \Deriv_{(y,y)} \sigma(0,\dd{C_2}{t}(0)) \nonumber \\
&=& 
-a \Deriv_{(y,y)} \sigma(0,\dd{C_1}{s}(0)) + 
b \Deriv_{(y,y)} \sigma(0,\dd{C_2}{t}(0)).
\label{eq:Derivative}
\end{eqnarray}
(We have used the chain rule and the derivative of the relation that
$\sigma(C_1(s),C_1(s))\equiv 0$.)
Transversality of the intersection of $C_1$ and $C_2$ at $0$ ensures
that the image of this differential is 
$$\Deriv_{(y,y)}\sigma(0\oplus T\Sigma_{y});$$
so transversality
of the section corresponding to $\QuotMap(y\times
p_2\times...\times p_k)$ 
at its zero $y$ ensures that the image of the
differential surjective onto the fiber of $E$ over $y$; i.e. the
canonical section is transverse.
The sign is correct, as one can see by inspecting
Equation~\eqref{eq:Derivative}.
\end{proof}

\begin{remark}
With the help of the above results, we can describe explicitly the
invariant of the tube:
$$\SWTube\colon H^*(\CritMan)\longrightarrow H^*(\Jac)\subset
\Alg(\Xtrunc),$$ 
which we do now for completeness.  Let
$\Lambda=\Wedge^*H_1(\Sigma)\subset \Alg(Y)$.  According to the proof
of Lemma~\ref{lemma:Commutator}, $\SWTube$ is a homomorphism of
$\Lambda$-modules; so, since $\Alg(Y)=\Lambda[U]$ surjects onto
$H^*(\CritMan)$, the invariant is determined by $\SWTube(U^i)$, as $i$
ranges over the non-negative integers. Since the Poincar\'e dual of
$\Sym^k(\Sigma)\subset T^{2g}$ (which is the image of
$\UnparModFlow(\CritMan,\Jac)$ under $\Restrict_{\Jac}$, according to
Theorem~\ref{thm:MOY}), is $$\frac{(\sum_{i=1}^g
\mu(A_i)\mu(B_i))^k}{k!},$$ it follows that
$$\SWTube(U^\ell)=\frac{(\sum_{i=1}^g A_i\cm
B_i)^{k+\ell}}{(k+\ell)!}.$$ We will not use this formula,
however. The results we prove in this paper require only the general
properties of $\SWTube$ which follow from its definition, together
with Proposition~\ref{prop:MovePoint}.
\end{remark}

According to Lemma~\ref{lemma:Commutator}, if $\gamma\subset Y$ is a
curve which is null-homologous in $\Xtrunc$, then it annihilates the
relative invariants, in the sense that
$$\SW_{(\Xtrunc,\CritMan)}(a\otimes \mu(\gamma)\cm\omega)=0.$$ If
sufficiently many curves in $Y$ become null-homologous in $\Xtrunc$,
then any class of sufficiently high degree in $\Alg(Y)$ annihilates
the invariant, as follows.

\begin{prop}
\label{prop:Algebra} Fix natural numbers $k,\ell$ with $\ell\geq k$.
Let $I$ denote the ideal generated by $\mu(A_1),\mu(A_2), ... ,
\mu(A_\ell)$ in $H^*(\Sym^k(\Sigma))$.  Then every element of
$H^*(\Sym^k(\Sigma))$ of degree greater than $k$ lies in $I$.
\end{prop}

\begin{proof}
The vector space $H^*(\Sym^k(\Sigma))$ is generated by homogeneous
elements of the form $$ U^{a}\cm \prod_{q=1}^{b}(A_{i_q}\cm
B_{i_q})\cm \prod_{r=1}^{c} A_{i_{b+r}}\cm \prod_{s=1}^{d}
B_{i_{b+c+s}}, $$ where $\{i_1,...,i_{b+c+d}\}$ is a subset of
$\{1,...,g\}$, and $a$, $b$, $c$, $d$ are integers with $a + b + c + d \leq k$.

Clearly, it suffices to prove the proposition for homogeneous 
generators of degree $k+1$.
Modulo $I$, such an element is equivalent to the element
$$\prod_{p=1}^{a}(U-A_p\cm B_p) \cm
  \prod_{q=1}^{b}(A_{i_q}\cm B_{i_q}-A_{a+q}\cm B_{a+q}) \cm 
  \prod_{r=1}^{c} A_{i_{b+r}}\cm 
  \prod_{s=1}^{d} B_{i_{b+c+s}}.$$
Indeed, in light of the fact that $a+b\leq k-d \leq \ell-d$, we can
arrange (after possibly simultaneously permuting the indices of the 
$\{A_i\}_{i=1}^g$ and $\{B_i\}_{i=1}^g$) that for each $s=1,...d$, 
$a+b<i_{b+c+s}$. Moreover, the original homogeneous element would
automatically lie in $I$ unless we had that $a+b<k\leq \ell<i_{b+r}$
for all $r=1,...,c$. Put together, must consider elements of the above
form which satisfy the constraint that $a+b < i_{j}$ for all
$j> b$.
If the degree of such an element is $k+1$, it must vanish in
$H^*(\Sym^k(\Sigma))$. 

This vanishing can be seen geometrically: $U$ is Poincar\'e dual to
the subset (identified with $\Sym^{k-1}(\Sigma)$) of $\Sym^k(\Sigma)$
where one point is constrained to lie in a specified point on
$\Sigma$: $A_i$ (resp. $B_i$) is Poincar\'e dual to the cycle where
one point is constrained to lie on $A_i$ (resp. $B_i$). Thus, (if one
chooses the point representing $U$ to be $A_i\cap B_i$), then
$U-A_i\cm B_i$ is Poincar\'e dual to the locus where two distinct
points are constrained; one is to lie on $A_i$, the other on $B_i$.
Similarly, the manifold Poincar\'e dual to $A_i\cm B_i-A_j \cm B_j$
gives a constraint on two distinct points in the symmetric
power. Finally, the remaining $A_{i_{b+r}}$ and $B_{i_{b+c+s}}$ give
additional, disjoint constraints (these are disjoint, if one chooses
that representing curves to be disjoint from the $A_i$ and $B_i$ for
$i=1,...,a+b$, which can be arranged since $a+b<i_{b+r}$ for all
$r\geq 1$).  Thus, since the total degree of the expression considered
is $k+1$, we have put constraints on $k+1$ distinct points, forcing the
intersection to be empty.
\end{proof}

\begin{cor}
\label{cor:VanishingOfRelInv}
Suppose that $\CritMan=\Sym^k(\Sigma)$, and
let $\ell\geq k$ be an integer so that there is a symplectic basis
$\{A_i,B_i\}_{i=1}^g$ for $H_1(\Sigma)$ so that $i_*(A_i)=0$ in
$H_1(\Xtrunc;\R)$ for $i=1,...,\ell$.  Then, for each $b\in\Alg(\Sigma)$ of
degree $d(b)>k$, and each $a\in\Alg(\Xtrunc)$, $\omega\in
H^*(\CritMan)$, we have 
$$\SW_{(\Xtrunc,\CritMan)}(a \otimes b \cm \omega)\equiv 0.$$
\end{cor}

\begin{proof}
By Proposition~\ref{prop:Algebra}, $b$ lies in the ideal generated by
$\mu(A_1),...,\mu(A_{\ell})$. Now the proposition follows from
Lemma~\ref{lemma:Commutator}.  
\end{proof}

We now have the promised proof of Proposition~\ref{prop:Vanishing}.

\vskip0.3cm
\noindent{\bf Proof of Proposition~\ref{prop:Vanishing}.}
Recall that we have constructed $\SW_{(\Xtrunc,\CritMan)}$ so that
$$\SWirr_{\spinc}(a \cm i_*(b))=\SW_{(\Xtrunc,\CritMan)}(a \cm i_*(b)
\otimes 1).$$
By Lemma~\ref{lemma:Commutator} and Proposition~\ref{prop:MovePoint},
we can write
$$\SW_{(\Xtrunc,\CritMan)}(a \cm i_*(b) \otimes 1)=
\SW_{(\Xtrunc,\CritMan)}(a\otimes b) + \SW_{(\Xtrunc,\Jac)}(a \cm c).$$
for some $c\in \Alg(\Sigma)$.  Note that $c$ lies in the ideal
generated by $H_1(\Sigma)$, as it can be expressed in terms of
Seiberg-Witten invariants of the tube, which take values in
$H^*(\Jac)\cong \Wedge^*(H_1(\Sigma))\subset
H_1(\Sigma)\cm\Alg(\Sigma)$.  By
Corollary~\ref{cor:VanishingOfRelInv}, the first term vanishes (using
the homological hypothesis of the inclusion of $\Sigma$ in $X$). The
remaining term is identified with an absolute invariant, according to
Proposition~\ref{prop:SmallSpinc}.
\hfill $\Box$

The proof of Proposition~\ref{prop:GeneralVanishing}, follows from the
same argument as Proposition~\ref{prop:Vanishing}; only in that case,
one must use the following (much simpler) analogue of
Corollary~\ref{cor:VanishingOfRelInv}.

\begin{lemma}
\label{lemma:GeneralVanishingOfRelInv}
Suppose that $\CritMan=\Sym^k(\Sigma)$. Then, for each $b\in\Alg(\Sigma)$ of
degree $d(b)>2k$, and each $a\in\Alg(\Xtrunc)$, $\omega\in
H^*(\CritMan)$, we have 
$$\SW_{(\Xtrunc,\CritMan)}(a \otimes b \cm \omega)\equiv 0.$$
\end{lemma}

\begin{proof}
This follows immediately from the fact that $\dim\CritMan=2k$.
\end{proof}
\section{The moduli spaces over $N$}
\label{sec:ModNSig}

The purpose of this section is to give the results about the
neighborhood of $\Sigma$ which were used in
Section~\ref{sec:ProductFormula}. Most of these results are
applications of~\cite{MOY} and \cite{SympThom}.  We assume for the
duration of this section that the $\SpinC$ structure over $\NSig$
satisfies $\langle c_1(\spinc),[\Sigma]\rangle \not\equiv n
\pmod{2n}$. We return to the excluded cases in
Section~\ref{sec:Perturbation}.

Over $\NSig$, endowed with a cylindrical-end metric and a certain
torsion connection on $T\NSig$ , the Seiberg-Witten equations admit a
complex interpretation analogous to the complex interpretation of the
equations over a K\"ahler manifold (see Section~5 of~\cite{SympThom}
for an explicit description of this connection, and especially
Proposition~5.6 where the complex interpretation is proved). The
Seiberg-Witten equations over $\NSig$ can be written as equations for
a connection $A$ over $E$,
$\alpha\oplus\beta\in(\Omega^{0,0}\oplus\Omega^{0,1})(\NSig,E)$:
\begin{eqnarray}
2 \Lambda F_A - \Lambda F_{\Canon{\NSig}}
&=& \frac{i}{2}(|\alpha|^2-|\beta|^2) 
\label{eq:VortexCurv}
\\
\Tr F_A^{0,2} &=& {\overline \alpha}\otimes \beta 
\label{eq:VortexInteg}\\
\DBar_A \alpha + \DBar_A^*\beta &=&0,
\label{eq:VortexHarm}
\end{eqnarray}
where $\Lambda$ denotes projection onto the $(1,1)$ form of the metric.
As noted in~\cite{SympThom},
for finite energy solutions,  decay estimates justify the usual
integration-by-parts which shows that one of $\alpha$ or $\beta$ must
vanish identically; i.e. the solutions over $\NSig$ correspond to
vortices over $\NSig$.

When $\beta\equiv 0$, then $A$ induces an integrable $\DBar$-operator
on $E$, $\DBar_A$, with respect to which $\alpha$ is
holomorphic. Moreover, by the usual exponential decay results,
together with the understanding of the solutions over $Y$
(Theorem~\ref{thm:MOY}), $(A,\alpha)$ exponentially approaches the
pull-back of a vortex solution over $\Sigma$. According to~\cite{MOY},
the underlying holomorphic data extends to the ruled surface $\RSurf$
obtained by attaching a copy of $\Sigma$ (denoted $\Sigma_+$)
to $\NSig$ ``at infinity.'' We state the results here for convenience.

\begin{defn}
Let $\Phi\in\Sections(\NSig,W^+)$, $\Psi\in\Sections(Y,W)$ be a pair
of spinors, and $\delta>0$ be some real number. Then, $\Psi$ is said
to {\em $\delta$-decay to $\Psi$} if for each $k\geq 0$,
$$\lim_{t\goesto\infty}\sup_{\{t\}\times Y}e^{\delta
t}|\nabla^{(k)}\Psi-\nabla^{(k)}\pi^*(\Psi)|=0,$$
where $\nabla^{(k)}$ denotes the $k$-fold covariant derivative.
More generally, $\Phi$ is said to {\em
decay to $\Psi$} if there is some $\delta>0$ so that $\Psi$
$\delta$-decays to $\Psi$. A similar notion can be defined for objects
other than spinors, such as connections, differential forms, etc.
\end{defn}

\begin{defn}
Given a line bundle $E$ over $Z$, a {\em holomorphic pair $(A,\alpha)$
in $E$} is a pair consisting of a $\DBar$-operator $\DBar_A$ over
$E$, and a section $\alpha$ of $E$, so that $F_A^{0,2}=0$, and
$\DBar_A\alpha=0$.
\end{defn}

\begin{theorem}
\label{thm:Extend}
Let $(A,\alpha)$ be a holomorphic pair on $\NSig$ which decays to a
the pull-back of a holomorphic pair $(A_0,\alpha_0)$ over
$\Sigma$. Then, there is a naturally associated line bundle ${\widehat
E}$ over $\RSurf$ and holomorphic pair $({\widehat A}, {\widehat
\alpha})$ in ${\widehat E}$, so that $(\DBar_{\widehat
A},{\widehat\alpha})|_{\RSurf-\Sigma_+}\cong (\DBar_A,\widehat
\alpha)$ and $(\DBar_{\widehat A},{\widehat\alpha})|_{\Sigma_+}\cong
(\DBar_{A_0},\alpha_0)$.
\end{theorem}

The above theorem is essentially a restatement of Theorems~7.7
of~\cite{MOY}, where it is stated for the cylinder, thought of as
$\RSurf$ minus two copies of $\Sigma$, rather than the neighborhood of
$\Sigma$, thought of as $\RSurf$ minus one copy of $\Sigma$ (though
the proof is no different).  Analogous results for the anti-self-dual
equations were obtained by Guo~\cite{Guo}.

In a similar vein we have the
following result, which allows us to deal with solutions with
reducible boundary values. We state the result slightly differently
from the above, since we will apply it in other contexts later.

\begin{theorem}
\label{thm:ExtendRed}
Let $A$ be a connection on a line bundle $E$ over $\NSig$, with $F_A^{0,2}=0$
and $E|_{(0,\infty)\times Y}\cong \pi^*(E_0)$, so
curvature form $F_A$ decays to the pull-back of a closed two-form
$F_0$ over $\Sigma$ with 
$$\Big(\frac{i}{2\pi}\int_\Sigma F_0\Big)-\langle c_1(E_0),[\Sigma]\rangle \not\in n\Z.$$  Then,
there is an associated line bundle ${\widehat E}$ over $\RSurf$ and
complex structure $\DBar_{\widehat A}$ with $$({\widehat
E},\DBar_{\widehat A})|_{\RSurf-\Sigma_+}\cong (E,\DBar_A),$$ and
$\langle c_1({\widehat E}),[\Sigma_+]\rangle$ is the greatest integer
congruent to $\langle c_1(E),[\Sigma]\rangle$ moduli $n$ smaller than
$\frac{i}{2\pi}\int F_0$. Furthermore, there is a natural identification
$$\Ker\Dirac_A\cap L^2 \cong H^0(\RSurf,{\widehat E})\oplus
H^2(\RSurf,{\widehat E}),$$ and $$\CoKer\Dirac_A\cap L^2 \cong
H^1(\RSurf,{\widehat E}).$$
\end{theorem}

The above is proved in Proposition~9.2 (see Corollary~9.11 and
Theorem~10.6) of~\cite{MOY}.

These results allow us to rule out the existence of certain solutions.
Recall first the following standard fact about the cohomology of
$\RSurf$ (see for example~\cite{Hartshorne}):

\begin{prop}
\label{prop:RSurfCohomology}
Let $\RSurf$ denote the ruled surface over $\Sigma$, which is given as
the projectivization of $\C\oplus L$, $\Projectivize(\C\oplus L)$
(here, $L$ is some line bundle over $\Sigma$).  Let ${\widehat E}$ be
a line bundle over the ruled surface $\RSurf$ and let $E_0$ denote the
restriction of of ${\widehat E}$ to $\Projectivize(\C\oplus 0)\cong
\Sigma$ and let $\ell$ be the evaluation of $c_1({\widehat E})$ on a
fiber in the ruling.  Then, if $\ell\geq 0$, $$ H^0(\RSurf,{\widehat
E})
\cong
\sum_{j=0}^{\ell}H^0(\Sigma,\ESigma\otimes\NormBundle^{\otimes j});~H^1(\RSurf,{\widehat E})
\cong
\sum_{j=0}^{\ell}H^1(\Sigma,\ESigma\otimes\NormBundle^{\otimes j});~H^2(\RSurf,{\widehat E})=0;
$$
and
if $\ell<0$, 
$$
H^0(\RSurf,{\widehat E})=0;~H^1(\RSurf,{\widehat E}) \cong
\sum_{j=1}^{-\ell-1}H^0(\Sigma,\ESigma\otimes\NormBundle^{\otimes
-j});~H^2(\RSurf,{\widehat E}) \cong 
\sum_{j=1}^{-\ell-1}H^1(\Sigma,\ESigma\otimes\NormBundle^{\otimes -j}).
$$
In particular, if $\ell=-1$, then $H^*(\RSurf,{\widehat E})=0$.
\end{prop}

We can apply these results to the case where $\NSig$ is a neighborhood
of a surface of self-intersection number with $-\Sigma\cm\Sigma >2g-2$. Recall that according to Theorem~\ref{thm:MOY}, for each $\SpinC$
structure $\spinc$ on $\NSig$, there are at most two components to the moduli
space of the boundary, the reducible component $\Jac$ and the
irreducible component $\CritMan$.

\begin{prop}
\label{prop:UnobstructedOnNNoPert}
If $$\Big|\langle c_1(\spinc),[\Sigma]\rangle\Big|<n,$$ 
the moduli space $\ModSW_N(\Jac)$ contains only reducibles. Moreover, the
space of reducibles is smoothly identified with the Jacobian torus
$\Jac$ (i.e. the kernel and the cokernel of the Dirac operator coupled
to any reducible vanishes). Furthermore, $\ModSW_N(\CritMan)$ is
empty.
\end{prop}

\begin{proof}
We prove that both moduli spaces $\ModSWirr_N(\Jac)$ and
$\ModSW_N(\CritMan)$ are empty.
Suppose there were some finite energy
solution to the Seiberg-Witten equations in a $\SpinC$ structure
with $|\langle c_1(\spinc),[\Sigma]\rangle|<n$. We know that
the spinor lies entirely in one of the two summands in the splitting
of the spinor bundle $W^+\cong E\oplus (\Canon{\NSig}^{-1}\otimes E)$
(i.e. it is an $\alpha$- or a $\beta$-spinor, in the notation of
Equations~\eqref{eq:VortexCurv}-\eqref{eq:VortexHarm}). By conjugating
if necessary (which switches the two summands and sends the $\SpinC$
structure $\spinc$ to another one $J\spinc$ with
$c_1(J\spinc)=-c_1(\spinc)$), we can assume without loss of generality
that the solution is an $\alpha$-solution. 

According to Theorem~\ref{thm:Extend} (and
Theorem~\ref{thm:ExtendRed}, when the boundary value is reducible), we
can extend the data $(E,\DBar_A, \alpha)$ over the associated ruled
surface $\RSurf$, obtained by attaching the curve $\Sigma_+$ at
infinity.  The fact that ${\widehat E}$ is an extension of $E$ says
that
\begin{eqnarray*}
\langle c_1({\widehat E}), [\Sigma_-]\rangle
&=& \langle c_1(E), [\Sigma]\rangle \\
&=& \OneHalf \langle c_1(\spinc)+c_1(\Canon{\NSig}),[\Sigma]\rangle
\\
&=& g-1 + \frac{n+\langle c_1(\spinc),[\Sigma]\rangle }{2} 
\end{eqnarray*}
where $\Sigma_-$ is the curve in ${\widehat E}$ with self-intersection
number $-n$ (which is identified with $\Sigma\subset \NSig$). By our
hypothesis, then, $$g-1< \langle c_1({\widehat E}),[\Sigma_-]\rangle
<n+g-1.$$ On the other hand, Equation~\eqref{eq:VortexCurv} says that
$\frac{i}{2\pi} F_A$ converges to the pullback of a form over $\Sigma$
whose integral is $g-1$, so Theorem~\ref{thm:ExtendRed} guarantees
that the Chern number of restriction to the other section of the
ruling satisfies the bound $$-n+g-1< \langle c_1({\widehat
E}),[\Sigma_+]\rangle < g-1.$$ Now, since the Poincar\'e dual of a
fiber is $(\PD[\Sigma_+]-\PD[\Sigma_-])/n$, we see that the evaluation
of $c_1({\widehat E})$ on a fiber is given by $$\ell = \frac{\langle
c_1(\widehat E),[\Sigma_+]\rangle - \langle c_1(\widehat
E),[\Sigma_-]\rangle}{n} =-1.$$ According to
Proposition~\ref{prop:RSurfCohomology}, it follows that
${\widehat\alpha}$ (and hence also $\alpha$) must vanish identically,
contradicting the irreducibility hypothesis on $(A,\alpha)$.

The fact that the reducibles are smoothly cut out in this range
follows in an analogous manner, using Theorem~\ref{thm:ExtendRed} and Proposition~\ref{prop:RSurfCohomology}.
\end{proof}

\begin{remark}
Most of this result can be found in Proposition~2.5 of~\cite{SympThom}.
\end{remark}

The above vanishing result is special to the particular $\SpinC$
structures considered, as it used the fact that the Dolbeault
cohomology of certain line bundles over the ruled surface vanish. In
general, the moduli spaces over $\NSig$ typically do contain
irreducibles. To study the deformation theory around these
irreducibles, we use an infinitesimal version of
Theorem~\ref{thm:Extend}; but first, we pause for a brief discussion
of deformation theory for the Seiberg-Witten equations in general.

In general, on a four-manifold $\Xtrunc$ with a cylindrical end, the
deformation complex around a solution $(A,\Phi)$ whose boundary value is
smooth and irreducible, is given by $$
\Omega^0(\Xtrunc,i\R) 
\longrightarrow
\Omega^1(\Xtrunc,i\R) \oplus \Sections(\Xtrunc,W^+)
\longrightarrow
\Omega^+(\Xtrunc,i\R)\oplus \Sections(\Xtrunc,W^-).
$$ Here, terms in $\Omega^0(\Xtrunc, i\R)$ are required to lie in
$\Sobol{2}{\delta,k}$, the $\delta$-decaying Sobolev space with $k$
derivatives (here we can choose any $k\geq 3$);
i.e. functions for which $$(\|f\|_{\delta,k})^2=\int_{\Xtrunc} (|f|^2
+ |\nabla f|^2 +...+ |\nabla^{(k)}f|^2)e^{\delta\tau} <\infty,$$ where
$\tau$ is a smooth function on $\Xtrunc$ which agrees with the $t$
coordinate over the cylindrical end.  Terms in $\Omega^1(\Xtrunc,i\R)
\oplus
\Sections(\Xtrunc,W^+)$ are required to lie in $\Sobol{2}{\delta,k-1}$
extended by the tangent space to the moduli space at infinity
at $\Restrict[A,\Phi]$. Finally, terms
in $\Omega^+(\Xtrunc,i\R)\oplus \Sections(\Xtrunc,W^-)$ are required to lie
in $\Sobol{2}{k-2}$. The first map in the deformation complex is the
linearization of the gauge group action on $[A,\Phi]$ around the
identity, while the second is the linearization of the Seiberg Witten
equations around $[A,\Phi]$.  When the boundary value of $[A,\Phi]$ is
a smooth reducible, then the above specifies the deformation theory
for the moduli space based at infinity.
In either case, the moduli space of solutions
about $[A,\Phi]$ is transversally cut out by the Seiberg-Witten
equations on $\Xtrunc$ if the $H^2$ of the above complex vanishes.
(This discussion is modeled on the theory
developed in~\cite{MMR}.)

The space of divisors in a compact, complex surface $X$ admits a
deformation theory, defined as follows. Consider the pair
$(\DBar_{\widehat A},{\widehat\alpha})$ where ${\DBar_{\widehat A}}$ is an integrable
$\DBar$-operator, and ${\widehat\alpha}$ is ${\DBar_{\widehat A}}$-holomorphic;
i.e.
\begin{eqnarray*}
F_{\widehat A}^{0,2} &=& 0 \\
\DBar_{\widehat A}\alpha &=& 0.
\end{eqnarray*}
This has a deformation complex
$$\begin{CD}
\Omega^{0,0}@>>>\Omega^{0,1}\oplus \Omega^{0,0}(E) @>>>
\Omega^{0,2}\oplus \Omega^{0,1}(E) @>>> \Omega^{0,2},
\end{CD}$$
whose cohomology groups are identified with the cohomology groups of
the quotient sheaf ${\mathcal E}/{\widehat\alpha}$, obtained from the
short exact sequence of sheaves: $$\begin{CD} 0@>>> {\mathcal O}_X
@>{\widehat \alpha}>> {\mathcal E} @>>> {\mathcal E}/{\widehat\alpha}
@>>>0.
\end{CD}
$$

\begin{theorem}
\label{thm:ExtendDeformations}
Let $(A,\alpha)$ correspond to a solution to the Seiberg-Witten
equations over $\NSig$, with irreducible boundary values. Then, the
cohomology groups of the deformation complex of the Seiberg-Witten
deformation complex are naturally isomorphic to the cohomology groups
deformation complex of the divisor $[\DBar_{\widehat
A},{\widehat\alpha}]$ in the line bundle ${\widehat E}$ over $\RSurf$,
provided by Theorem~\ref{thm:Extend}. When $(A,\alpha)$ has a
reducible boundary value, then $H^2$ of the Seiberg-Witten deformation
complex is identified with $H^1(\RSurf,{\mathcal
E}/{\widehat\alpha})$, while the tangent space of the based moduli
space is identified with $\C\oplus H^0(\RSurf,{\mathcal
E}/{\widehat\alpha})$. 
\end{theorem}

\begin{proof}
This follows exactly as in Theorem~9.14 (for irreducible boundary
values) and Theorem~10.12 (for reducible boundary values)
of~\cite{MOY}. The key observation at this point is to note that
$$\Lambda\Del\DBar+|\alpha|^2\colon \Sobol{2}{\delta,k}\longrightarrow
\Sobol{2}{\delta,k-2}$$
is an isomorphism, which allows one to ``unroll'' parts of the
Seiberg-Witten deformation complex to identify it with the deformation
theory of divisors in $\NSig$. As in \cite{MOY} (see
also~\cite{KMStructureTheorem}), we can identify $\Lambda\Del\DBar$
with the operator over the cylindrical end with
$$-e^{-2\lambda t}\DDt e^{2\lambda t}\DDt + \Delta_Y,$$
where $\lambda=\frac{\pi n}{\Vol(\Sigma)}$. According to the theory
of~\cite{LockhartMcOwen}, the operator 
$$\Lambda\Del\DBar\colon \Sobol{2}{k,\delta} \longrightarrow
\Sobol{2}{k-2,\delta}$$
is Fredholm for all weights $0<\delta<4\lambda$. In particular, it has
the same index for all small $0<\delta$ as it has on the weight
$\delta=2\lambda$, where it can be connected via Fredholm operators to
the manifestly self-adjoint operator $$d^{*_\lambda}d\colon
\Sobol{2}{k,\lambda}\longrightarrow
\Sobol{2}{k-2,\lambda},$$
where
$d^{*_\lambda}$ denotes the formal $\lambda$-weighted adjoint of $d$. 
It follows from the homotopy invariance of the index that
$\Lambda\Del\DBar+|\alpha|^2$ has index zero on
$\Sobol{2}{k,\delta}$. 
From the maximum principle, it has no kernel, so it induces an
isomorphism as claimed, identifying the deformation theory of the
Seiberg-Witten equations with the deformation theory of divisors in
$\NSig$. Passing to the ruled surface then follows from Corollary~9.4
of~\cite{MOY}. 
\end{proof}

\begin{prop}
\label{prop:NSigDefTheory}
Let $\NSig$ be a disk bundle over a surface $\Sigma$ with
$$\Sigma\cm\Sigma = -n < 2-2g,$$ 
endowed with a $\SpinC$ structure $\spinc$ with
$$n < |\langle c_1(\spinc),[\Sigma]\rangle| \leq
n+2g-2.$$ 
Let $$e=\frac{n+2g-2-|\langle
c_1(\spinc),[\Sigma]\rangle|}{2}.$$ Then, the expected dimensions of the
moduli spaces over $N$ and $\Xtrunc$ are given by:
\begin{eqnarray}
\edim \ModSW_N(\Jac)&=& 2\mult+1 \\
\edim \ModSW_N(\CritMan)&=& 2\mult.
\end{eqnarray}
Moreover, $\ModIrr_N(\Jac)$, $\ModSW_N(\CritMan)$, are transversally
cut out by the Seiberg-Witten equations (in particular, they are
manifolds of the expected dimension). Furthermore, the boundary map
$$\Restrict\colon\ModSW_N(\CritMan)\longrightarrow\CritMan$$
is an orientation-preserving diffeomorphism onto its image.
\end{prop}

\begin{proof}
The deformation theory around a solution $[A,\alpha]\in\ModSW_N(\CritMan)$ is
identified the deformation theory around a corresponding divisor in
the line bundle ${\widehat E}$ with $\langle
c_1({\widehat E}),[\Sigma_{\pm}]\rangle=e$; i.e. with
divisors in a line bundle which (topologically) pulls back from
$\Sigma$. According to Proposition~\ref{prop:RSurfCohomology},
all such divisors actually pull back
from the base $\Sigma$; and indeed, the deformation theory corresponds to
deformation theory of degree $e$ divisors in the base $\Sigma$, which
is unobstructed. Thus,  $\ModSW_N(\CritMan)$ is a manifold of
real dimension $2\mult$, transversally cut out by the Seiberg-Witten
equations. 

The above transversality applies to $\ModSW_N(\Jac)$ as well, except
that the expected dimension is greater by one, as we saw in
Theorem~\ref{thm:ExtendDeformations}.

This identification of deformation theories of $\ModSW_N(\CritMan)$
proves that $\Restrict$ is an orientation-preserving local
diffeomorphism onto its image in 
$\Sym^{\mult}(\Sigma)\cong \CritMan$. In fact, it
is injective, as follows. As we saw, any two solutions with the same
boundary values actually vanish over the same disks (with the same
multiplicities). By the usual analysis of the vortex equations, any
two such solutions must differ by a complex gauge transformation;
i.e. a function $u$ which satisfies $$\Lambda \DBar\Del u + |\alpha|^2
(e^{2u}-1)=0,$$ where $u$ is a function which decays on the
cylinder. By the maximum principle, such a function must vanish
identically.
\end{proof}

Having analyzed the moduli spaces over neighborhoods of $\Sigma$, we
close with a some general results concerning the rest of the moduli
spaces of the complement of $\Sigma$.

\begin{prop}
\label{prop:XtruncDefTheory}
Let $\Xtrunc$ be as in Proposition~\ref{prop:DimensionCounting}. Then,
letting $\edim \ModSW(X)=d$, we have 
\begin{eqnarray}
\edim \ModSW_{\Xtrunc}(\Jac) &=& d+2g-2\mult-2 \\
\edim \ModSW_{\Xtrunc}(\CritMan) &=& d.
\end{eqnarray}
Moreover, 
$\ModSW_{\Xtrunc}(\Jac)$, and $\ModSW_{\Xtrunc}(\CritMan)$ are
transversally cut out by the Seiberg-Witten equations (in particular,
they are manifolds of the expected dimension).
\end{prop}

\begin{proof}
By a standard excision argument,
we have
$$\edim \ModSW_{\Xtrunc}(\Jac) + \edim \ModSW_N(\Jac) - 2g + 1 = \edim
\ModSW_X(\spinc)=d,$$
which calculates $\edim\ModSW_{\Xtrunc}(\Jac)$, given
Proposition~\ref{prop:NSigDefTheory}. Similarly, we have
$$\edim \ModSW_{\Xtrunc}(\CritMan) + \edim \ModSW_N(\CritMan) - 2\mult
= d,$$
which gives us $\edim\ModSW_{\Xtrunc}(\CritMan)$.

Smoothness of $\ModSW_{\Xtrunc}(\Jac)$ and $\ModSW_{\Xtrunc}(\CritMan)$
follows from adapting methods of~\cite{MMR}.
\end{proof}

\section{Perturbations when $e=g-1$}
\label{sec:Perturbation}

In our earlier discussion, we had to exclude one $\SpinC$ structure
over $Y$. In this section, we introduce a perturbation of the
equations which allows us to handle this case. We begin by adapting
results of Section~\ref{sec:ModY} to this perturbed equation, and
then, we give a discussion which is parallel to that of
Section~\ref{sec:ModNSig}. The perturbations used here are analogues
of Taubes' perturbations in the symplectic
category~\cite{TaubesSympI}, \cite{TaubesSympII}; see also~\cite{KMO} for a related discussion.

Recall that the Seiberg-Witten equations over $Y$ are obtained as the critical 
points of the Chern-Simons-Dirac functional
$\CSD$ defined over the configuration space
$$\SWConfig(Y,W)=\Conns(W)\oplus \Sections(Y,W)/\Map(Y,S^1),$$
where $\Conns(Y,W)$ denotes the space of 
connection in the spinor bundle $W$ which are compatible with some
fixed connection $\TConn$ on $TY$.  
The
functional is defined by $$\CSD(B,\Psi)=\int_Y (B-B_0)\wedge \Tr(F_B+
F_{B_0}) -
\int_Y \langle \Psi, \Dirac_B \Psi
\rangle, $$ where $B_0\in\Conns(Y,W)$ is some reference connection,
$B-B_0\in\Omega^1(Y;i\R)$ denotes the difference $1$-form, and $\Tr$
denotes the trace of the corresponding connection on $W$. Its
Euler-Lagrange equations (the three-dimensional Seiberg-Witten
equations) are
\begin{eqnarray}
*\Tr (F_B) - i\tau(\Psi) = 0 \label{eq:CurvEqThree}\\
\Dirac_B \Psi = 0, \label{eq:HarmSpinThree} 
\end{eqnarray}
where
$$\tau\colon \Gamma(Y,W)\rightarrow \Omega^{1}(Y;\R) $$
is adjoint 
to Clifford multiplication, in the sense that for all $\gamma\in 
\Omega^{1}(Y;\R)$, $\Psi\in \Gamma(Y,W)$, we have 
\begin{equation}
\label{eq:DefOfTau}
\frac{1}{2}\langle i\gamma\cm \Psi,\Psi\rangle_W
=-\langle \gamma, \tau(\Psi)\rangle_{\Lambda^1}.
\end{equation}
Moreover, its upward gradient flow equations are the usual
Seiberg-Witten equations on the four-manifold $\R\times Y$.

When $Y$ is a circle-bundle over a Riemann surface
with Euler number $-n$ satisfying
$$n>2g-2,$$ recall that these
equations are inconvenient in the $\SpinC$ structure when $e=g-1$
(in the notation of Section~\ref{sec:ModY}). We will find it
useful to consider a perturbed functional 
$$\CSD_u\colon \SWConfig(Y,\spinct)\longrightarrow
\R,$$
where $u\in\R$, 
given by $$\CSD_u(B,\Psi)= \CSD(B,\Psi)+ u \int_Y i\eta\wedge(\Tr F_B
- \Tr F_{B_0}),$$ where $\eta$ is the connection form for $Y$ over
$\Sigma$, and the reference connection $B_0$ satisfies
$\Tr(F_{B_0})\equiv 0$ (i.e. $B_0\in\Jac$).
These give rise to perturbed Seiberg-Witten equations of the form
\begin{eqnarray}
*\Tr (F_B) - i\tau(\Psi) + i u (*d\eta) = 0 \label{eq:CurvEqPert}\\
\Dirac_B \Psi = 0, \label{eq:HarmSpinPert}
\end{eqnarray}
whose moduli space of solutions is denoted $\ModSWThree_u(Y)$. The
gradient flow equations of the perturbed functional are solutions to
the Seiberg-Witten equations on $\R\times Y$, perturbed by the
self-dual component of $iu (d\eta)$, which can be collected into
moduli spaces, denoted
$\ModFlow_u(\CritMan_1,\CritMan_2)$, or their unparameterized versions
$$\UnparModFlow_u(\CritMan_1,\CritMan_2)=
\ModFlow_u(\CritMan_1,\CritMan_2)/\R.$$
We have the following analogue of Theorem~\ref{thm:MOY}.

\begin{theorem}
\label{thm:MOYPert}
Let $Y$ be a circle-bundle over a Riemann surface with genus $g>0$ and
Euler number $-n<2-2g$. Let $\spinct$ be the $\SpinC$ structure
corresponding to $g-1\in \Zmod{n}\subset H^2(Y;\Z)$.  For all $u$ with
$0<u<2$, the moduli space contains two components, a reducible one
$\Jac$, identified with the Jacobian torus $H^1(\Sigma;\R/\Z)$, and a
smooth irreducible component $\CritMan$ diffeomorphic to
$\Sym^{g-1}(\Sigma)$. Both of these components are non-degenerate in
the sense of Morse-Bott. There is an inequality
$\CSD_{u}(\Jac)>\CSD_{u}(\CritMan)$, so the space
$\ModFlow_{u}(\Jac,\CritMan)$ is empty.  The space
$\UnparModFlow_{u}(\CritMan,\Jac)$ is smooth of expected
dimension $2g-2$; indeed it is diffeomorphic to $\Sym^{g-1}(\Sigma)$.
\end{theorem}

\begin{proof}
Most of this is a straightforward adaptation of~\cite{MOY}. 

We begin with the identification of the moduli spaces over $Y$.  As
in~\cite{MOY}, the equations over $Y$ reduce to vortex equations over
$\Sigma$. More specifically, the components of the moduli spaces
$\ModSWThree_Y(\spinct)$ correspond to line bundles $E_0$ over
$\Sigma$ with the property that
$$\pi^*(E_0\oplus\Canon{}^{-1}\otimes E_0)\cong W,$$ the spinor bundle
of $\spinct$ (here $\Canon{}$ denotes the canonical line bundle
over $\Sigma$). The vortex equations are are equations for
$B\in\Conns(E_0)$, $\alpha\oplus \beta\in \Sections(\Sigma, E_0 \oplus
\Canon{}^{-1}\otimes E_0)$, 
which, in the case at hand, take the form 
\begin{eqnarray}
2 F_B - F_\Canon{} + i u (d\eta)  &=&
i(|\alpha|^2-|\beta|^2)(*1) 
\label{eq:VortCurvEqPert}\\
\DBar_B \alpha + \DBar_B^*\beta &=& 0 
\label{eq:VortHarSpinPert}\\
\alpha\otimes \beta &=& 0.
\end{eqnarray}
Thus, one of $\alpha$ or $\beta$ must vanish. In fact, in our case,
$$\deg E_0\equiv g-1\pmod{n}.$$
In fact, if $$\deg E_0\neq g-1,$$ then the
solution space to these equations ($0<u<2$) is empty. More
specifically, letting $\deg E_0=g-1+n\ell$, we see that
when $\beta\not\equiv 0$, then by integrating 
Equation~\eqref{eq:VortCurvEqPert}  over $\Sigma$
against $i/2\pi$, we get $$2(g-1+n\ell)-(2g-2)+u\deg Y
=2n\ell-un \geq 0,$$ which forces $\ell\geq 1$ (since
$u>0$). 
Since in this case $\deg(E_0)>g-1$, $H^1(\Sigma,E_0)=0$, so $\beta$ must vanish. 
If, on the other other hand, it is $\alpha\not\equiv
0$, then we obtain in the same manner that $$2n\ell-un
\leq 0,$$ which forces $\ell\leq 0$ (since $u<2$). Since $\alpha$ 
represents a class in $H^0(\Sigma,E_0)$, it follows that $\ell=0$.

So, all irreducibles correspond to $\alpha$-vortices in the line
bundle $E_0$ with $\deg E_0=g-1$.  The identification of this space
of vortices with the symmetric product follows from~\cite{Bradlow}
(see also~\cite{JaffeTaubes}).

Non-degeneracy of the irreducible manifold $\CritMan$ follows exactly
as in~\cite{MOY}. To see non-degeneracy of $\Jac$, we appeal to
results of Section 5.8  of~\cite{MOY}. Consider
the Dirac operator on the $\SpinC$ structure with spinors
$W=E\otimes (\C\oplus\Canon{}^{-1})$ with connection induced
from a connection $B\in\Conns(E)$ whose curvature pulls up from
$\Sigma$.
It is shown in Proposition~5.8.4 of~\cite{MOY} that 
this Dirac operator
admits no harmonic spinors unless the holonomy around a fiber circle
in $Y$ is trivial. In fact this holonomy is trivial when the following integral is congruent to $g-1$ modulo $n\Z$:
$$\frac{i}{4\pi^2}\int_Y F_B\wedge\eta = 
g-1-
\frac{u\deg(Y)}{2},$$ (we have used here
Equation~\eqref{eq:CurvEqPert}). Since $0<u<2$, this holonomy is
non-trivial, so the reducibles admit no harmonic spinors, i.e. $\Jac$
is smoothly cut out by the equations.

We now perform the Chern-Simons calculations (see
the proof of Proposition~5.23 of~\cite{MOY}). Suppose
$[(B_1,\Psi_1)]\in\CritMan$, and $[(B_0,0)]\in \Jac$. Then, we have 
\begin{eqnarray*}
2\deg B_0 - \deg \Canon{} + u\deg(Y)&=&0; \\
2\deg B_1 - \deg \Canon{} &=& 0,
\end{eqnarray*}
where by $\deg B$, we mean the integral $\frac{i}{4\pi^2}\int_Y
F_B\wedge \eta$, which when $B$ is induced from a line bundle over
$\Sigma$, agrees with the degree of that line bundle.
So,
\begin{eqnarray*}
\CSD_u(B_1) &=& \int_Y (B_1-B_0)\wedge (2F_{B_1} + 2F_{B_0} - 2
F_{\Canon{}}) 
+ u \int i \eta \wedge (2 F_{B_1}- 2 F_{B_0}) \\ 
&=&\frac{8\pi^2}{\deg Y}(\deg B_1 - \deg B_0)(\deg B_1 + \deg B_0
- \deg \Canon{}) \\
&&+ u \int i \eta \wedge (2 F_{B_1}- 2 F_{B_0}) \\ 
&=& 2\pi^2 u^2 \deg Y,
\end{eqnarray*}
which is negative; while
$\CSD_u(B_0)=0$. 

The smoothness of the space of flows, and its identification with the
symmetric product, follows exactly as in the unperturbed case (see
Section~\ref{sec:ModY}).
\end{proof}

We now turn to the neighborhood of $\Sigma$. We use a perturbation
over $\NSig$ which is compatible with the above perturbation over
$Y$. Specifically, let $$f\colon \NSig \longrightarrow \R$$ be a smooth
function which is identically zero on the complement of the cylinder
$[0,\infty)\times Y\subset
\NSig$, and identically one on the subcylinder $[1,\infty)\times Y$.
We consider the Seiberg-Witten equations perturbed by the self-dual
part of $i u f (d\eta)$.
Note that this perturbing
two-form is $i u
\lambda f$ times the $(1,1)$ form of the standard cylindrical-end
metric on $\NSig$ (see~\cite{SympThom}), where $$\lambda=-\frac{2\pi
\deg Y}{\Vol(\Sigma)}$$ (here, $\Vol(\Sigma)$ denotes the volume of
$\Sigma$). Similarly, we can extend the perturbation over $Y$ to a
self-dual two-form perturbation of the equations over $\Xtrunc$ (and,
consequently, $X(T)$ to all $T>2$).  Denote the corresponding
moduli spaces by $\ModSW_{N,u}(\Jac)$, $\ModSW_{N,u}(\CritMan)$,
$\ModSW_{\Xtrunc,u}(\Jac)$, $\ModSW_{\Xtrunc,u}(\CritMan)$,
and $\ModSW_{X(T),u}$. Strictly speaking, we still have to show that
these perturbed moduli spaces $\ModSW_{X(T),u}(\spinc)$ can be used to
calculate the Seiberg-Witten invariant in either chamber. This is
clear because we can always choose a compactly-supported 
perturbing two-form $\eta_0$ whose integral against $\omega_g$
dominates the integral of $\omega_g$ against $iuf(d\eta)^+$. The key
point is that the latter integral is finite, since $\omega_g$ decays
exponentially (see~\cite{APSI}).

We now have the following analogue of
Proposition~\ref{prop:UnobstructedOnN}

\begin{prop}
\label{prop:UnobstructedOnNPert}
Suppose $\langle c_1(\spinc),[\Sigma]\rangle=n$, and let $u$ be a
real number with $0<u<2$. Then the perturbed moduli space
$\ModSW_{N,u}(\Jac)$ contains only reducibles. Moreover, the space of
reducibles is smoothly identified with the Jacobian torus $\Jac$
(i.e. the kernel and the cokernel of the Dirac operator coupled to any
reducible vanishes). Furthermore, $\ModSW_{N,u}(\CritMan)$ is empty.
\end{prop}

\begin{proof}
We begin by proving $\ModSW_{N,u}(\CritMan)$ is empty.  Note that
$\CritMan$ consists entirely of $\alpha$-solutions, hence so must any
section in $\ModSW_{N,u}(\CritMan)$. Thus, a solution
$(A,\alpha)\in\ModSW_{N,u}(\CritMan)$ induces a non-zero element in
$H^0({\widehat E})$ with
\begin{eqnarray*}
\langle c_1({\widehat E}),[\Sigma_-]\rangle = n+g-1 &{\text{and}}&
\langle c_1({\widehat E}),[\Sigma_+]\rangle = g-1.
\end{eqnarray*}
But $H^*(\RSurf,{\widehat E})\equiv 0$, according to
Proposition~\ref{prop:RSurfCohomology}. The same argument, now
appealing to Theorem~\ref{thm:ExtendRed}, shows that
$\ModSWIrr_{N,u}(\Jac)$ is empty, and that $\Jac$ is smooth.
\end{proof}

\begin{prop}
\label{prop:DimensionCountingPert}
Suppose that 
$$\langle c_1(\spinc),[\Sigma]\rangle =-n,$$ 
and let $u$
be any real number with $0<u<2$. Then according to
Theorem~\ref{thm:MOYPert}, $\ModSWThree_u(\spinc|_Y)$ has two
components, $\Jac$ and $\CritMan$, where $\CritMan$ is diffeomorphic
to $\Sym^{g-1}(\Sigma)$.  Furthermore, the expected dimensions of the
moduli spaces over $N$ and $\Xtrunc$ are given by:
\begin{eqnarray}
\edim \ModSW_{N,u}(\Jac)&=& 2g-1 \\
\edim \ModSW_{N,u}(\CritMan)&=& 2g-2 \\
\edim \ModSW_{\Xtrunc,u}(\Jac) &=& 2d \\
\edim \ModSW_{\Xtrunc,u}(\CritMan) &=& 2d.
\end{eqnarray}
Moreover, $\ModIrr_{N,u}(\Jac)$, $\ModSW_{N,u}(\CritMan)$,
$\ModSW_{\Xtrunc,u}(\Jac)$, and $\ModSW_{\Xtrunc,u}(\CritMan)$ are
transversally cut out by the Seiberg-Witten equations (in particular,
they are manifolds of the expected dimension).
Furthermore, the boundary map
$$\Restrict\colon\ModSW_{N,u}(\CritMan)\longrightarrow\CritMan$$
is an orientation-preserving diffeomorphism onto its image.
\end{prop}

\begin{proof}
The proofs of Propositions~\ref{prop:NSigDefTheory} and \ref{prop:XtruncDefTheory}
apply directly in this perturbed context.
\end{proof}

\section{Cohomology}
\label{sec:Cohomology}

The Seiberg-Witten invariant is obtained from pairings of certain canonical
cohomology classes on the Seiberg-Witten moduli space. These cohomology
classes are inherited from the configuration spaces in which the
moduli spaces live. In this section, we recall the definitions of these
classes and discuss natural geometric representatives for
them. (See Chapter~5 of~\cite{DonKron} for the corresponding
discussion of cohomology relevant to Donaldson invariants.)

Let $X$ be a Riemannian four-manifold with a $\SpinC$ structure $\spinc$
specified by the pair of Hermitian $\C^2$ bundles $W^+$ and $W^-$, and
the Clifford action $$\rho\colon TX\otimes W^+\longrightarrow W^-.$$
The Seiberg-Witten pre-configuration space  is the space
$$\SWPreConfig(W^+)=
\Conns(W^+)\times \Sections(X;W^+)\cong \Omega^1(W;\R)\times
\Sections(X;W^+),$$
where $\Conns(W^+)$ denotes the space of connections compatible with
some fixed connection $\TConn$ on $TX$, and the isomorphism above is
induced by comparing any connection $A$ against some fixed connection
$A_0$.  The irreducible pre-configuration space $\SWPreConfigIrr(W^+)$ is
the subset of $\SWPreConfig(W^+)$ consisting of pairs $(A,\Phi)$, where
$\Phi\not\equiv 0$. Now, $\SWPreConfigIrr(W^+)$ is weakly contractable, and
the space $\Map(X;S^1)$ acts freely on it, so the {\em irreducible
configuration space}, which is
$$\SWConfigIrr(W^+) = \SWPreConfigIrr(W^+)/\Map(X;S^1)$$
is weakly homotopy equivalent to the classifying space of
$\Map(X;S^1)$. Now,
$$\Map(X;S^1)\sim \Map(X;S^1)_{e}\times \pi_0(\Map(X;S^1))
\sim S^1\times H^1(X;\Z);$$
so 
$$\SWConfig(W^+)\sim \CP{\infty}\times \frac{H^1(X;\R)}{H^1(X;\Z)},$$
and
$$H^*(\SWConfig(W^+);\Z)\cong \Z[U]\otimes \Wedge^* H^1(X;\Z),$$
where $U$ is a generator with grading two. 
More invariantly, we define
$$\Alg(X)=\Z[H_0(X;\Z)]\otimes \Wedge^*H_1(X;\Z),$$
graded by declaring $H_0(X;\Z)$ to have grading two and $H_1(X;\Z)$ to
have grading one. Then, we have seen that
$$H^*(\SWConfig(W^+);\Z)\cong \Alg(X).$$

We describe two functorial mechanisms for constructing 
generators in $H^*(\SWConfigIrr(W^+);\Z)$. 
Over the space $X\times\SWConfigIrr(W^+)$, there is a universal
line bundle ${\mathcal L}=X\times S^1 \times
\SWPreConfigIrr(W^+)/\Map(X,S^1)$, where the action is defined by
$$u(x,\zeta,A,\Phi)=(u, u(x)\zeta, A + u^{-1}du, u\Phi).$$
Using this class we can define a ``$\mu$-map''
$$\mu\colon (H_0\oplus H_1)(X;\Z)\longrightarrow H^*(\SWConfigIrr(W^+)),$$
which sends a homology class of degree $i$ to a cohomology class of
degree $2-i$, by
$$\mu(x)=c_1({\mathcal L})/x;$$
i.e. $\mu(x)$ is the cohomology class on $\SWConfigIrr(W^+)$ with the
property that for any homology class $c\in H_*(\SWConfigIrr(W^+))$,
$$\langle \mu(x), c\rangle = \langle c_1({\mathcal L}), x\times
c\rangle.$$

We describe another convenient mechanism for constructing
one-dimensional cohomology in $\SWPreConfig(W^+)$ as follows. A closed
curve $x\colon S^1\longrightarrow X$ induces a map $$\Hol_x\colon
\SWConfig(W^+)\longrightarrow S^1$$ which is defined to be the
holonomy of the connection $A$ around the curve $x$. The pull-back of
the volume form $d\theta$ of $S^1$ by this map gives rise to a
one-dimensional cohomology class $\Hol_x^*(d\theta)$ associated to
$x$, which we call the {\em holonomy class around $x$}.

\begin{prop}
\label{prop:MuCircle}
The cohomology groups of the configuration space $\SWConfig(W^+)$ are
generated by the image of the $\mu$-map. Moreover, given $x\in
H_1(X;\Z)$, $\mu(x)$ is the holonomy class around $x$, 
$\Hol_x^*(d\theta)|_{\SWConfigIrr(W^+)}$.
\end{prop}

\begin{remark} Note that $\Hol_x^*(d\theta)$ is naturally defined over
the entire configuration space
$$\SWConfig(W^+)=\SWPreConfig(W^+)/\Map(X,S^1)\sim \frac{H^1(X;\R)}{H^1(X;\Z)}.$$
\end{remark}

\begin{proof}
We begin by proving the second claim.
Note that ${\mathcal L}$ comes with a tautological connection along
the $X$ factor, with the property that for any path $\beta\colon
S^1\longrightarrow X$ and connection $A\in
\SWPreConfig(W^+)$, 
\begin{equation}
\label{eq:UniversalConnection}
\Hol_{\beta\times A}({\mathcal L})=\Hol_{\beta}(A).
\end{equation}
$\SWPreConfig(W^+)$
Now, fix a path in $X$, $$\beta\colon S^1\longrightarrow X.$$ We need
to show that for all paths in the configuration space $$\alpha\colon
S^1\longrightarrow \SWConfig(W^+)(\spinc),$$ we have that $$\langle
c_1({\mathcal L}), {\alpha\times\beta}\rangle = \deg
(\Hol_{\beta}\circ \alpha\colon S^1\longrightarrow S^1).$$

This follows from the fact that for a line bundle $L$ over the the
torus $S^1\times S^1$, the first Chern number is the degree of the map
from $S^1\times S^1$ defined by $$x\mapsto \Hol_{x\times S^1}L$$ (a
map which makes sense only after one puts a connection on $L$, but the
degree is independent of this connection, so we left it out of the
notation), together with the universal property of
Equation~\eqref{eq:UniversalConnection}. Thus, we have identified
$\mu$ on any one-dimensional homology class. 

The rest of the proposition is established, once we see that for a
point $x\in X$, $\mu[x]$ generates $H^2$ of the configuration
space. But this follows easily from the fact that $\Map(X,S^1)_e$ acts
freely on the space of irreducible configurations.
\end{proof}

With this concrete understanding of the $\mu$-classes, we turn to a
discussion of submanifold representatives for them.

Given a point $x\in X$, let ${\mathcal L}_x$ denote the line bundle
associated to the base fibration of $X$; i.e. it is the restriction of
the universal line bundle ${\mathcal L}$ to the slice
$\SWConfigIrr(W^+)\cong \{x\}\times \SWConfigIrr(W^+)\subset X\times
\SWConfigIrr(W^+)$. 
Given a point in the fiber $\Psi(x)\in W^+_x$, we can construct a 
canonical section over $\SWConfig(X,W^+)$ by
$$[A,\Phi]\mapsto [A,\Phi,\langle \Phi(x),\Psi(x)\rangle].$$
The zero set of this section in
$\SWConfig(X,W^+)$ is a codimension-two submanifold representing
$\mu[x]$. The restriction of this section to a moduli space
$\ModSW_X(\spinc)\subset\SWConfig(X,\spinc)$ is not, in general,
transverse. However, by mollifying the construction appropriately, we
can find a section which is generic over the moduli space, and hence
obtain a divisor $\Divisors(x)$ representing $\mu[x]$, as follows.

\begin{defn}
Fix a ball $B\subset X$ around $x$ and a non-vanishing section
$\Psi$ of $W^+|_B$. Given a self-dual two-form $\lambda$ which is
compactly supported over $B$, the {\em $\lambda$-mollified section} is
the section of ${\mathcal L}_x$ defined by
$$[A,\Phi]\mapsto [A,\Phi,\int_B \langle
\lambda\cm\Phi,\Psi\rangle].$$
\end{defn}

\begin{lemma}
There are $L^2$ sections $\lambda$ compactly supported in $B$ so that
the 
$\lambda$-mollified section of ${\mathcal L}_x$, restricted to the
moduli space $\ModSW_{X}$, vanishes transversally. 
\end{lemma}

\begin{proof}
Fix a compactly-supported cut-off function $\beta$ in $B$, and
consider the section
$$[A,\Phi]\times \lambda 
\mapsto [A,\Phi,\int_B \langle
\lambda\cm\Phi,\Psi\rangle\beta]$$
of $\pi_2^*({\mathcal L}_x)$, thought of as a line bundle over $\Omega^+(B)\times
\ModSW(X,\spinc)$, giving $\Omega^+(B)$ the $L^2$ topology. This
transversality follows from the fact that, for any
$[A,\Phi]\in\ModSW_X$, as we vary $\lambda$, the
integral $\int_B \langle
\lambda\cm\Phi,\Psi\rangle\beta$ can take on any complex value. This,
in turn, follows from the unique continuation theorem for elliptic
differential operators, which guarantees that the section $\Phi$
cannot vanish identically over $B$.
\end{proof}

\begin{remark}
In effect, the above lemma tells us how to construct a divisor
representative $\Divisors(x)$ for $\mu[x]$ when $[x]\in H_0(X)$; this
divisor is represented by the zero-set of the $\lambda$-mollified
section of ${\mathcal L}_x$. Finding codimension-one representatives
for $\mu[\gamma]$, where $\gamma\in H_1(X)$ is even easier: one need
only find a regular value $\theta$ for the map
$$\Hol_\gamma\colon\ModSW_X(\spinc)\rightarrow S^1.$$ Then,
$\Hol_\gamma^{-1}(\theta)$ is the submanifold $\Divisors(\gamma)$
representing $\mu[\gamma]$ over $\ModSW_X(\spinc)$.
\end{remark}

\bibliographystyle{plain}
\bibliography{biblio}

\end{document}